\documentclass[11pt,a4paper]{article}
\usepackage{setspace}
\usepackage[utf8]{inputenc}
\usepackage[T1]{fontenc}
\usepackage{geometry}
\usepackage{fancybox}
\usepackage{enumerate}
\usepackage{epsfig}
\usepackage{amsmath} 
\usepackage{amsfonts}
\usepackage{amssymb}
\usepackage{amsthm}
\usepackage{marvosym}
\usepackage{color}
\usepackage{import}
\usepackage{natbib}
\usepackage{bbold}
\usepackage{adjustbox}
\usepackage{booktabs, multirow}
\usepackage{authblk}
\usepackage{todonotes}
\usepackage{hyperref}

\usepackage{algorithm} 
\usepackage{algpseudocode}
\usepackage{supertabular}

\algdef{SE}[DOWHILE]{Do}{doWhile}{\algorithmicdo}[1]{\algorithmicwhile\ #1}%

\newtheorem{definition}{Definition}
\newtheorem{proposition}{Proposition}
\newtheorem{remark}{Remark}

\def \N {\mathbb{N}}
\def \Z {\mathbb{Z}}

\definecolor{couleurOB}{rgb}{1.0 , 0.0 , 0.0}
\definecolor{couleurHC}{rgb}{0.0 , 0.5 , 0.0}
\definecolor{couleurNC}{rgb}{0.0 , 0.0 , 1.0}
\definecolor{couleurAL}{rgb}{1.0 , 0.6 , 0.0}
\definecolor{couleurMO}{rgb}{0.5 , 0.0 , 0.5}
\definecolor{couleurDC}{rgb}{0.0 , 0.5 , 0.5}
\definecolor{darkcyan}{rgb}{0.0, 0.55, 0.55}

\newcommand{\pbname}{JOBPRSP-D}

\newcommand{\tour}{route}

\newcommand{\capa}{B}
\newcommand{\pool}{\mathcal{K}_{pool}}
\newcommand{\algo}{CGH}

\newcommand{\rev}[1]{{\textcolor{black}{#1}}}
\newcommand{\reven}[1]{{\textcolor{black}{#1}}}
\newcommand{\ie}{\emph{i.e.}\ }

\newcounter{constraint}
\newcommand{\constraintlabel}[1]{\refstepcounter{constraint}\label{#1}(\theconstraint)}

\DeclareMathAlphabet\mathbfcal{OMS}{cmsy}{b}{n}

\begin{document}

\title{A bin packing formulation for the joint batching, routing and sequencing problem}

\author[1]{Olivier Briant}
\author[1]{Hadrien Cambazard}
\author[2,3]{Diego Cattaruzza}
\author[1]{Nicolas Catusse}
\author[4]{Anne-Laure Ladier}
\author[2]{Maxime Ogier}

\affil[1]{\small CNRS, Grenoble INP, G-SCOP, University of Grenoble Alpes, Grenoble 38000, France}
\affil[2]{\small Univ. Lille, CNRS, Inria, Centrale Lille, UMR 9189 CRIStAL, F-59000 Lille, France}
\affil[3]{\small Department of Mathematics, Computer Science and Physics, University~of~Udine, Via delle Scienze 206, 33100 Udine, Italy}
\affil[4]{\small INSA Lyon, Université Lumière Lyon 2, Université Claude Bernard Lyon 1, Université Jean Monnet Saint-Etienne, DISP UR4570, 69621 Villeurbanne, France}
\affil[ ]{\small \texttt{\{olivier.briant, hadrien.cambazard, nicolas.catusse\}@grenoble-inp.fr}}
\affil[ ]{\small \texttt{diego.cattaruzza@uniud.it}}
\affil[ ]{\texttt{\small anne-laure.ladier@insa-lyon.fr}}
\affil[ ]{\small \texttt{maxime.ogier@centralelille.fr}}

\date{\today} 

\maketitle

\begin{abstract}
    Warehouses are 
    the scene of complex logistic problems integrating different decision layers. This paper addresses the Joint Order Batching, Picker Routing and Sequencing Problem with Deadlines (\pbname) in rectangular warehouses. To tackle the problem, an exponential linear programming formulation related to bin packing is proposed. It is solved with a column generation heuristic able to provide lower and upper bounds on the optimal value. We start by showing that the {\pbname} is related to a bin packing problem as opposed to a scheduling problem. We take advantage of this aspect to derive a number of valid inequalities that enhance the resolution of the master problem. The proposed algorithm is evaluated on publicly available data-sets. It is able to optimally solve instances with up to 18 orders in \reven{a} few minutes. It is also able to prove optimality or to provide high-quality lower bounds on larger instances with 100 orders as long as the picker's capacity is small enough. Finally, new upper bounds are found on all the instances except two of the benchmark considered. To the best of our knowledge this is the first article that provides optimality guarantee on large size instances for the {\pbname}; the results can therefore be used to assert the quality of heuristics proposed for the same problem.
    
\end{abstract}


\section{Introduction}

\begin{figure}[hbt]
\centering
\includegraphics[width=0.75\textwidth]{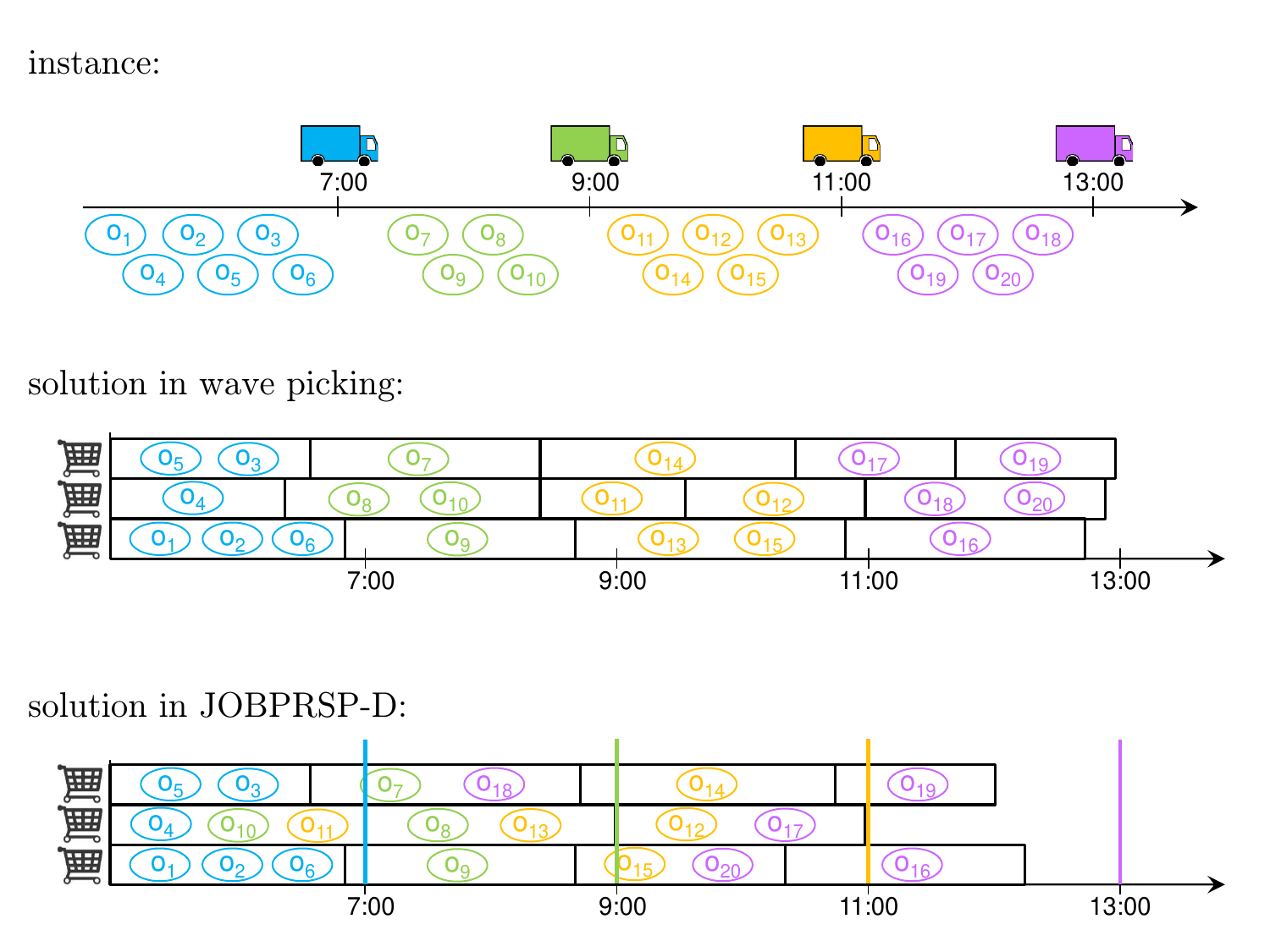}
 \label{fig:sequencing}
\caption{Interest of a sequencing problem when orders have deadlines.}
\end{figure}

The scene takes place in a tall and dim building, lights flickering with a buzzing noise. Pickers move across the warehouse, collecting items, fulfilling \rev{customers'} orders, with a digital voice for companionship.
This paper takes a narrow viewpoint restricted to the optimization problem of efficient warehouse operations management. To increase the efficiency of the operations, several orders might be batched, \ie put together to be collected by a picker in a single route. Batching is possible as long as orders grouped together respect capacity constraints, namely they fit in a trolley. Since customer satisfaction is at the core of competitiveness in e-commerce, each order is required to be prepared and shipped before a strict deadline that usually corresponds to the departure of a truck from the warehouse to deliver orders. Pickers must therefore collect their batches in an appropriate sequence in such a way that all orders meet their deadlines. Planning the entire process is referred to as the Joint Order Batching, Picker Routing and Sequencing Problem with Deadlines (\pbname). In the present work, we consider optimal routing policies in a rectangular warehouse with several blocks (or cross-aisles). Such optimal routes can be efficiently computed with dynamic programming following the work of \cite{Ratliff1983,Pansart18}. 
An example of the problem addressed in this work is provided in Figure~\ref{fig:sequencing} where each deadline is associated with a color.
In this instance, there are several sets of orders each associated with a deadline.
In a classical wave picking approach, the solution would be to batch together orders with the first deadline before batching together orders with the second deadline and so on.
Hence in the proposed example with three pickers, each picker sequentially retrieve\reven{s} batches with orders of the same deadline.
On the contrary, in the \pbname, orders of different deadlines can be mixed in the same batch, as long as for each deadline orders have been collected on time. The resulting solution has \rev{the} smallest total picking time as depicted in the figure.

Order picking in warehouses is part of a complex logistic process involving a wide range of problems such as storage assignment, order batching, picker routing and picker scheduling (\cite{DeKoster2007}). Relatively efficient methods have been designed for single problems such as the picker routing problem. However\rev{,} integrated problems are receiving attention since significant gains can be achieved by taking into account the interrelations of the planning problems (\cite{VanGils2018}). The resulting combinatorial problems embed many NP-Hard sub-problems with distinct decision levels which make them very challenging to solve in practice. Moreover, practical instances tend to be of large size which complicates the resolution even more. Consequently, heuristic approaches have been favored so far and, as mentioned in the review of \cite{VanGils2018} (page 6), {\em the use of mathematical programming as a research method to integrate different order picking planning problems is limited}. 

\cite{BriantCCCLO20} propose a column generation heuristic able to produce lower and upper bounds for the Joint Order Batching and Picker Routing Problem (JOBPRP). In this work, we show that the methodology proposed by \cite{BriantCCCLO20} can be extended to tackle the sequencing question as well, \ie the additional assignment and sequencing of batches to each picker in order to meet the deadlines of the orders. This aspect is currently seen as a scheduling problem in the literature (\cite{VanGils2019, Haouassi2022}). We take a different viewpoint and propose a \emph{bin packing} formulation. As a result, we benefit from families of valid inequalities identified for packing problems and derived from Dual-Feasible Functions (see \cite{Clautiaux2010}). A high-quality lower bound for the integrated problem can therefore be derived by adding such cuts in the formulation proposed in \cite{BriantCCCLO20}. The bin packing viewpoint also provides new insights to design efficient heuristics. Finally, a key cutting plane that was used for solving the pricing of~\cite{BriantCCCLO20} can be considerably strengthened by using timing considerations related to the sequencing aspect of the problem.

The contribution\reven{s} of the present paper are: 
\begin{itemize}
    \item an extended mathematical formulation based on bin packing considerations is proposed for the {\pbname}. As already stated, the {\pbname} is currently seen as a scheduling problem in the literature (\cite{VanGils2019, Haouassi2022}). This bin packing viewpoint can also be of inspiration for future researchers to develop ad-hoc methods for the {\pbname};
    \item solving the extended formulation yields proven optimal solutions on all the small benchmark instances for the {\pbname};
    \item families of valid inequalities for the {\pbname} are derived from Dual-Feasible Functions;
    \item tour constraints proposed in \cite{BriantCCCLO20} are strengthened based on the particular structure of the benchmark of instances for the {\pbname};
    \item the column generation based heuristic proposed is able to provide high-quality upper and lower bounds for subsets of the benchmark instances;
    \item by extending the methodology proposed in \cite{BriantCCCLO20}, we show its generality and potential to efficiently tackle (and provide lower bounds) for more integrated problems related to the JOBPRP.
\end{itemize}


The paper is organized as follows. Section~\ref{sec:pbspec} \reven{specifies the problem and notations}. A literature review is presented in Section~\ref{sec:litrev}. Formulations and valid inequalities are presented in  Section~\ref{sec:bpform} -- this is the key section on the paper; a reader who is pressed for time can focus on Section \ref{sec:bpform}. Section~\ref{sec:algo} \reven{then details} each component of the proposed methodology as well as the overall algorithm design. Experimental results are reported and discussed in Section~\ref{se:expes}. Finally Section~\ref{sec:conclusion} \reven{draws} some conclusions. 

\section{Problem definition and notation\reven{s}}
\label{sec:pbspec}

The warehouse layout is modeled as a directed graph  $\mathcal G=(\mathcal{V}, \mathcal{A})$ where $\mathcal{V}$ contains two depots $s$ and $s'$ and two types of vertices: locations $\mathcal{V}_L$ and intersections $\mathcal{V}_I$. So, overall $\mathcal{V} = \mathcal{V}_L \cup \mathcal{V}_I \cup \{s, s'\}$. Each location in set $\mathcal{V}_L$ contains one or more products to be picked, whereas the intersections in $\mathcal{V}_I$ are used to encode the warehouse structure. Additionally, $s$ denotes the depot where picking {\tour}s start and $s'$ is the depot where {\tour}s end \reven{to drop off} picked orders.
Moreover, a travelling time $t_{ij}$ is associated with each arc $ (i, j)\in \mathcal{A}$. 

The set of items is denoted by $\mathcal{I}$. 
An item $i \in \mathcal{I}$ is characterized by its location $v_i \in \mathcal{V}_{L}$. 
Note that here we consider an item has a unique location.

A set $\mathcal{O}$ of orders to collect is given. Each order $o \in \mathcal{O}$ is composed of a set $\mathcal{I}_o \subseteq \mathcal{I}$ of items to pick.
Let us denote by $\capa_o = |\mathcal{I}_o|$ the size of order $o\in\mathcal{O}$. 

An order $o\in\mathcal{O}$ is associated with a time $t_o^\text{pick}$ that represents the time needed to pick the items composing order $o$. It is proportional to the size $B_o$ of order $o$.
We indicate as $\mathcal{V}_o\subseteq\mathcal V_{L}$ the set of locations that have to be visited to retrieve order $o\in\mathcal O$.
Moreover, an order $o \in \mathcal O$ is associated with a deadline $\bar{d}_o > 0$ that represents the latest  instant of time to have the order ready at node $s'$. 

The set of deadlines is denoted as $\mathcal{D}$, and in practice the number of different deadlines in $\mathcal{D}$ is much lower than the number of orders since deadlines correspond to the departure of a vehicle from the warehouse to deliver orders.
Hence many orders share the same deadline.

A batch is a set of orders that are collected together in a single {\tour} by a picker. Each batch is associated to a {\tour} $k$ in the warehouse starting from $s$, terminating in $s'$ and visiting a subset of locations in $\mathcal V_L$. The orders collected along route $k$ are denoted as $\mathcal{O}_k \subseteq \mathcal{O}$. The deadline $\bar{d}_k$ of a route $k$ is the minimum deadline among the deadlines of all its orders, that is $\bar{d}_k = \min_{o \in \mathcal{O}_k} \bar{d}_o$.

The time $t_k$ to perform a route $k$ is the sum of three terms: the picking times of all order to be collected in the route, the travel times to go from each location to the next one, and, finally, a setup time $t^\text{setup}$. The setup time is supposed constant and represents the time needed by the picker to get the trolley ready.

The routes are performed by a set $\mathcal P$ of pickers. Each picker pushes a trolley of capacity $\capa$. A route $k$ must respect the trolley's capacity, \ie must satisfy $\sum_{o \in \mathcal{O}_k} \capa_o \leq \capa$.

A solution of the {\pbname} is a sequence of {\tour}s $\mathcal{K}_p$ for each picker $p\in\mathcal{P}$ such that:

\begin{itemize}
    \item all orders are picked up;
    \item each {\tour} respects the capacity constraint of the trolley;
    \item in each sequence $\mathcal{K}_p = \{r^1_p, r^2_p,\ldots, r^{n_p}_p\}$, {\tour} $r^{i+1}_p$ starts after {\tour} $r^i_p$ ends;
    \item each {\tour} arrives at the final depot no later than the minimum deadline of the collected orders.
\end{itemize}

The objective of the {\pbname} is to minimize the sum of the times to perform the determined {\tour}s.
Indeed, since the number of available pickers and orders to be collected is an input, minimizing the total time improves the productivity of the system. Moreover, \reven{since} deadlines \reven{are} considered here as hard constraints, there is no need to consider tardiness aspect\reven{s} in the objective function.

In the following, we will indicate by $\mathcal K$ the set of all possible feasible {\tour}s. Moreover, given two {\tour}s $k$ and $k'\in\mathcal K$, we may write $k\subseteq k'$ instead of $\mathcal{O}_k \subseteq \mathcal{O}_{k'}$ for ease of notation. Similarly, given an order $o\in\mathcal{O}$ and a route $k\in\mathcal K$, we may write $o\in k$ instead of $o\in\mathcal{O}_k$.





\section{Literature review}
\label{sec:litrev}

The literature on problems related to Order Batching and Picker Routing is large and an exhaustive review is \rev{out of scope} of the present work. The interested reader is referred to the literature reviews of \cite{DeKoster2007}, \cite{VanGils2018}, \cite{cergibozan2019}, \cite{Vanheusden2022}.

We thus limit the review presented in this section to works directly related to this paper.
First, we review works on the JOBPRP proposing algorithms that are exact or that produce both lower and upper bounds on the studied problem. 
Then, we review works that consider the sequencing of pickers.

\subsection{Solving the JOBPRP exactly or producing lower and upper bounds}

Some exact algorithms have been proposed for the JOBPRP but rarely extended to more complex integrated problems. 
Algorithmic techniques relevant to the present work have also been investigated on the JOBPRP and are mostly based on exponential Integer Programming (IP) formulations inspired by vehicle routing. 

In particular, to the best of our knowledge, \cite{Gademann2005} are the first to design a Branch-and-Price algorithm where variables are related to batches and their corresponding routes in the warehouse. The proposed approach remains limited to the case of a single block warehouse where all orders have the same size. The pricing algorithm is a dedicated Branch-and-Bound algorithm that uses the dynamic program of \cite{Ratliff1983} to compute the optimal route of a batch. \cite{Muter2015} generalize the pricing problem of \cite{Gademann2005} to deal with orders of different sizes, and also consider heuristic routing policies. \cite{Valle2017} propose a Branch-and-Cut approach to solve the JOBPRP in a \rev{two-block} warehouse and an optimal routing policy. The optimal routing is based on sub-tour elimination constraints typical of the Traveling Salesman Problem formulation by Dantzig. The algorithm proposed \reven{can} optimally solve instances \reven{up to} 20 orders. The dynamic programming algorithm of \cite{Ratliff1983} is then generalized to any number of cross-aisles and used by \cite{BriantCCCLO20} to design another column generation technique. The pricing is solved with integer programming and the addition of cutting planes based on the dynamic program for computing optimal routes. 
Recently, \cite{Hessler2022} propose to extend the state space procedure of \cite{Ratliff1983} to formulate the pricing problem as a shortest path problem with side constraints, namely the grouping requirement of orders and the trolley capacity. It also naturally handles scattered storage \ie an item may be stored at more than one pick location. The pricing problem can then be efficiently solved with an IP model. The proposed algorithm remains limited to a \rev{one-block} warehouse.
Finally, \cite{Wahlen2023} propose a Branch-and-Price-and-Cut algorithm to solve the JOBPRP.
The pricing problem is formulated as a Shortest Path Problem with Resource Constraints (SPPRC) on a linear graph where the nodes represent the customer orders in a given sequence, and two arcs between two consecutive nodes represent the inclusion or exclusion of the order in the batch.
The SPPRC is solved by a labelling algorithm with strong completion bounds; the results on instances with \rev{one-block} warehouses outperform \reven{those} of \cite{Muter2015}.

\subsection{Problems with picker sequencing considerations}

The sequencing of pickers, which is at the heart of the present paper, is closely related to the wave-based picking problem. Picking is often performed in the so-called wave-based manner where a group of orders (a wave) is released at the same time to be picked. The goal is to minimize the makespan of the wave which requires to assign batches to pickers. This is a natural extension of the JOBPRP where the exact sequence of each picker is not needed since the makespan is only dependent on the assignment of batches to the pickers and not on their sequencing. Wave picking problems with makespan minimisation are addressed by \cite{Ardjmand2018, Ardjmand2020, Muter2022}.


Sequencing really matters when \reven{considering} due dates or deadlines for orders. 
A number of variants of this integrated batching, routing and sequencing problem are addressed in the literature; they vary on the objective and the\reven{ir consideration of} due dates as either hard or soft constraints. 
Table~\ref{table:literature} provides an overview of the related literature.
Columns \texttt{time} and \texttt{tardiness} indicate if the objective function minimize\reven{s} the time and the tardiness respectively.
Column \texttt{Due date constraints} indicates if due dates are considered as hard or soft: note that a hard due date corresponds to a deadline.
Finally, column \texttt{Solution method} indicates the proposed solution methods.


\begin{table}[H]
\begin{adjustbox}{center}
  \footnotesize
    \begin{tabular}{lcccl}
    \toprule 
    \multicolumn{1}{c}{\multirow{2}[2]{*}{\bf Reference}} & \multicolumn{2}{c}{\bf Objective} & \multicolumn{1}{c}{\bf Due date} & \multirow{2}[2]{*}{\bf Solution method} \\ \cmidrule{2-3}
          & \bf time  & \bf tardiness & \bf constraints &  \\
    \midrule
    \cite{Tsai2008} & \checkmark & \checkmark & soft  & genetic algorithm \\
    \midrule
    \multirow{2}[2]{*}{\cite{Henn2013}} & \multirow{2}[2]{*}{} & \multirow{2}[2]{*}{\checkmark} & \multirow{2}[2]{*}{soft} & iterated local search \\
          &       &       &       & tabu search \\
    \midrule
    \multirow{2}[2]{*}{\cite{Henn2015}} & \multirow{2}[2]{*}{} & \multirow{2}[2]{*}{\checkmark} & \multirow{2}[2]{*}{soft} & variable neighborhood descent \\
          &       &       &       & variable neighborhood search \\
    \midrule
    \multirow{2}[2]{*}{\cite{Chen2015}} & \multirow{2}[2]{*}{} & \multirow{2}[2]{*}{\checkmark} & \multirow{2}[2]{*}{soft} & genetic algorithm \\
          &       &       &       & ant colony \\
    \midrule
    \cite{Scholz2017} &       & \checkmark & soft  & variable neighborhood descent \\
    \midrule
    \cite{Menendez2017} &       & \checkmark & soft  & variable neighborhood search \\
    \midrule
    \cite{VanGils2019} & \checkmark &       & hard  & iterated local search \\
    \midrule
    \cite{Haouassi2022} & \checkmark &       & hard  & heuristic + constraint programming \\
    \midrule
    \cite{DHaen2022} & \checkmark & \checkmark & soft  & large neighborhood search \\
    \bottomrule
    \end{tabular}%
    \end{adjustbox}
    \caption{Studies of variants of the integrated batching, routing and sequencing problem.}
    \label{table:literature}
\end{table}%

Table~\ref{table:literature} shows that several works focus on the sum of the delays to the due dates whereas others consider them as deadlines\rev{, \emph{i.e.},} as hard constraints, then focus on the minimization of the total time. 
Also, two papers consider the minimization of a weighted sum of tardiness and total time.
Additionally, note that two \rev{real-life} aspects of the problem are taken into account in recent works: dynamic arrival of orders (\cite{DHaen2022}) and scattered storage (\cite{Haouassi2022}).
All the algorithms proposed so far are based on metaheuristics such as genetic algorithm, iterated local search, tabu search, variable neighborhood search, ant colony, variable neighborhood descent, large neighbourhood search. 
Note that \cite{Haouassi2022} use a constraint programming approach to efficiently build a feasible schedule for each picker. 
\cite{Henn2015} mentions column generation as a possible resolution approach but turns to a heuristic. 

To the best of our knowledge, none of the algorithms proposed in the literature is exact or computes lower bounds for this integrated problem.
The algorithm proposed in this paper is meant to close this gap and contributes \reven{towards} assessing the numerous heuristics developed in this field. 

Moreover, \reven{Table~\ref{table:literature} also shows that considering hard deadlines for orders in the \pbname{} is a relevant set-up; it has been used by both \cite{VanGils2019} and \cite{Haouassi2022}, who conducted experiments over a shared and publicly available benchmark.}

\section{Formulation, relaxation and valid inequalities for the {\pbname}}
\label{sec:bpform}
This section first provide\reven{s} an extended formulation of the {\pbname} as a (multiple) bin packing problem. Let us first recall that given a {\tour} $k\in\mathcal K$, $t_k$ represents its total time, $\mathcal{O}_k \subseteq \mathcal{O}$ the batch of orders to be collected while performing {\tour} $k$ and $\bar{d}_k = \min_{o \in \mathcal{O}_k} \bar{d}_o$ its deadline. Additionally, to easily state the model, we denote as:
\begin{itemize}
\item $\mathcal{K}(o)=\{k\in\mathcal{K}:o \in\mathcal O_k\}$ the set of routes retrieving order $o$, for all orders $o \in  \mathcal{O}$;
\item  $\displaystyle\mathcal{K}(\bar{d})=\{k\in\mathcal{K}:\bar{d}_k\leq\bar{d}\}$ the set of routes with a deadline that is lower than or equal to $\bar{d}$, for all $\bar{d}\in\mathcal{D}$.
\end{itemize}

We use binary variables $\rho_{kp}$ that equal 1 if {\tour} $k\in\mathcal K$ is assigned to picker $p\in\mathcal P$, 0 otherwise. 
The {\pbname}  can be formulated as follows:


$$
(B(\mathcal{K}))\left\{
\begin{array}{rclllllllllllll}
\min&\displaystyle\sum_{k\in\mathcal{K}}t_k\;\sum_{p\in\mathcal{P}}\rho_{kp}&&&&\constraintlabel{objBP}&\\
\\
s.t.
&\displaystyle\sum_{k\in\mathcal{K}(o)}\; \sum_{p \in \mathcal{P}}\rho_{kp}\;\geq\; 1
&&\forall o\in\mathcal{O}&&\constraintlabel{cnst1BP}&\\
\\
&\displaystyle\sum_{k\in\mathcal{K}(\bar{d})}t_{k}\;\rho_{kp}\;\leq\; \bar{d}
&&\forall \bar{d}\in\mathcal{D}, \forall p\in\mathcal{P}&&\constraintlabel{cnst2BP}&\\
\\
&\rho_{kp}\;\in \{0,1\}
&&\forall k\in\mathcal{K}, p \in \mathcal{P}&&\constraintlabel{cnst3BP}&\\
\\
\end{array}
\right.
$$

The objective function~\eqref{objBP} minimizes the total time required by the pickers to retrieve all orders.
Constraints \eqref{cnst1BP} ensure that each order is collected. From a bin packing point of view, those constraints state that among the routes that compose the final packing there is at least one that contains a given order $o\in\mathcal{O}$.

Constraints~(\ref{cnst2BP}) impose the respect of the deadlines. In particular, for each deadline $\bar{d}\in\mathcal{D}$, they impose that the total time of {\tour}s in $\mathcal{K}(\bar{d})$ assigned to a picker do\reven{es} not exceed the value of $\bar{d}$. From a bin packing point of view, these constraints impose that the total size (\ie the picking time)  of the objects (\ie the {\tour}s) assigned to a bin (\ie a picker) do\reven{es} not exceed its capacity (\ie the deadline). Again here we emphasize that the sequence in which the picker performs the {\tour}s assigned to him/her is not important as long as their total time is lower than or equal to the considered deadline. This allows to look at the problem as a (variant of the multiple) bin packing problem.



Note  that a valid sequence for each picker is easily computed from an assignment of the routes to the pickers provided by $(B(\mathcal{K}))$ by sorting the routes by increasing deadlines.
Moreover, this bin packing formulation of the problem is possible because there are no release dates for the orders, \ie all orders can be collected from the beginning of the planning horizon.

Finally constraints~\eqref{cnst3BP} define the domain of the variables.

Note that this extended formulation contains an exponential number of variables with regards to the size of the instance, and solving directly such a formulation implies that all feasible routes have been generated beforehand.

\subsection{A relaxation of formulation $B(\mathcal{K})$}
\label{sec:colgenform}

Let us now consider the binary variables $\rho_k$ defined as:
$$
\displaystyle\rho_{k} = \sum_{p\in\mathcal{P}}\rho_{kp} \qquad \forall k\in\mathcal{K} \qquad \constraintlabel{cnst:rhokp}\\
$$

Note that $\rho_k$ equals 1 if {\tour} $k\in\mathcal K$ is selected to be part of the solution, 0 otherwise. 
By taking into account constraints~\eqref{cnst:rhokp} and aggregating constraints~\eqref{cnst2BP} in formulation $(B(\mathcal{K}))$, we obtain the following formulation:
$$
(M(\mathcal{K}))\left\{
\begin{array}{rclllllllllllll}
\min&\displaystyle\sum_{k\in\mathcal{K}}t_k\;\rho_{k}&&&&\constraintlabel{obj:relax}\\
\\
s.t.
&\displaystyle\sum_{k\in\mathcal{K}(o)}\;\rho_{k}\;\geq\; 1
&&\forall o\in\mathcal{O}&&\constraintlabel{cnst1:relax}&\\
\\
&\displaystyle\sum_{k\in\mathcal{K}(\bar{d})}t_{k}\;\rho_{k}\;\leq\; \bar{d}\times |\mathcal{P}|
&&\forall \bar{d}\in\mathcal{D}&&\constraintlabel{cnst2:relax}&\\
\\
&\rho_{k}\;\in \{0,1\}
&&\forall k\in\mathcal{K}&&\constraintlabel{vardef:relax}\\
\\
\end{array}
\right.
$$

The objective function~\eqref{obj:relax} minimizes the total time of the selected {\tour}s.
Constraints \eqref{cnst1:relax} ensure that each order is collected. Constraints \eqref{cnst2:relax} are obtained by summing constraints \eqref{cnst2BP} over all pickers $p \in \mathcal{P}$ and taking into consideration the relation defined in~\eqref{cnst:rhokp}.
Constraints~\eqref{vardef:relax} define the domain of the variables.
Note that $(M(\mathcal{K}))$ is the formulation proposed in \cite{BriantCCCLO20} for the JOBPRP with the addition \reven{of} constraints \eqref{cnst2:relax}.


\begin{remark}
Formulation $(M(\mathcal{K}))$ is a relaxation of $(B(\mathcal{K}))$.\\
This is a straightforward consequence of the definition of variables $\rho_{k}, k\in\mathcal{K}$ and the aggregation of constraints~\eqref{cnst2BP} to obtain constraints~\eqref{cnst2:relax}.
\end{remark}

\begin{proposition}
The linear relaxations of formulations $(B(\mathcal{K}))$ and $(M(\mathcal{K}))$ are equal.
\end{proposition}
\begin{proof}
Let us denote by $(LB(\mathcal{K}))$ and $(LM(\mathcal{K}))$ the linear relaxations of formulations $(B(\mathcal{K}))$ and $(M(\mathcal{K}))$ respectively: integrality requirements on the variables, \ie constraints \eqref{cnst3BP} and \eqref{vardef:relax}, are removed. 
From any fractional solution $\tilde\rho_{kp}$ of $(LB(\mathcal{K}))$, a feasible solution to $(LM(\mathcal{K}))$ can be built by setting $\tilde\rho_{k} = \sum_{p\in\mathcal{P}}\tilde\rho_{kp}$, that is using relation~(\ref{cnst:rhokp}). 
Reversely, from any solution $\tilde\rho_{k}$ of $(LM(\mathcal{K}))$, \reven{setting} $\tilde\rho_{kp} = \frac{1}{|P|}\tilde\rho_{k}$ for each picker $p \in P$ \reven{gives} a feasible solution of $(LB(\mathcal{K}))$. 
\end{proof}


Note that constraints \eqref{cnst2:relax} can be seen as enforcing in $(M(\mathcal{K}))$ the linear relaxation of the bin packing constraints~\eqref{cnst2BP} of $(B(\mathcal{K}))$. 
This bound is also known as the $L_1$ lower bound for bin packing. Stronger bounds (such as $L_2$) can be used to derive stronger inequalities following the framework of Dual-Feasible Functions.

\subsection{Valid inequalities for the {\pbname} from Dual-Feasible Functions}
\label{sec:ineqdff}

\vspace{1em}
In this section we derive valid inequalities for the {\pbname} by using Dual-Feasible Functions (DFF). The concept of DFF has been used to solve optimization problems involving knapsack constraints, like for example bin packing, cutting stock, scheduling or routing problems.
DFF are used either to compute lower bounds or to add valid inequalities in integer programs\footnote{A pedagogical video on DFF is available on this page: \url{https://youtu.be/5S5FOPW5mzA}}.
The interested reader is referred to \cite{Clautiaux2010} who propose a survey on this topic.
In the following, we recall some basic notions about DFF.

\begin{definition}
A function $f:\left[ 0,1\right] \rightarrow \left[0,1\right]$ is dual-feasible if for any finite set $S$ of real non-negative numbers, we have the relation $$\sum_{x \in S}{x} \leq 1 \;\; \Rightarrow \;\; \sum_{x \in S}{f(x)} \leq 1 .$$ 
\end{definition}

From any Dual-Feasible Function, the following proposition allows to derive valid inequalities for integer programs.

\begin{proposition}
If 
\begin{itemize}
    \item $f:\left[ 0,1\right] \rightarrow \left[0,1\right]$ is a DFF
    \item $S = \left\lbrace x \in \Z^n_+ ~ | ~ \sum_{j=1}^n{a_{ij} x_{j}} \leq b_i, \;\; \forall i=1, ..., m \right\rbrace$, with $b_i \geq a_{ij} \geq 0$, for all $i=1,\dots,m, j=1,\dots, n$ and $b_i>0$ for all $i = 1,\dots,m$
\end{itemize}
 then for any $i = 1,\dots,m$, $\sum_{j=1}^n{f(\frac{a_{ij}}{b_i}) x_j} \leq 1$ is a valid inequality for $S$.%
\label{prop:DFF}
\end{proposition}

\begin{proof}
See Section 2.3 in the paper of \cite{Clautiaux2010}.
\end{proof}

Numerous DFF for the bin packing problem have been proposed in the literature. We present in Section~\ref{sec:l2cuts} the so-called Martello and Toth's cuts and in Section~\ref{sec:FScuts} the so-called Fekete and Schepers's cuts.

\subsubsection{Martello and Toth's cuts}
\label{sec:l2cuts}

As explained above, each deadline in the definition of the {\pbname} can be seen as the size of a bin from a bin packing problem point of view. Thus, given a deadline $\bar{d} \in \mathcal{D}$, we can consider the $L_2$ bound for the related bin packing problem (\cite{Martello1990}). This provides a family of valid inequalities for formulation $(M(\mathcal{K}))$. 
Note that similar families of inequalities can be proposed for formulation $(B(\mathcal{K}))$.
Let us denote $\mathcal{Q}(\bar{d}) = \{q\in\mathbb{N} | q \geq 1, q \leq \left\lfloor \frac{1}{2}\bar{d} \right\rfloor\}$.
Given a deadline $\bar{d}\in\mathcal{D}$ and an integer $q\in\mathcal{Q}(\bar{d})$, we define the following sets:
\begin{itemize}
\item $\mathcal{K}_1(\bar{d}, q) = \left\lbrace k \in \mathcal{K}(\bar{d}) ~ | ~ t_k > \bar{d} - q \right\rbrace$ \reven{contains} all the {\tour}s in $\mathcal{K}(\bar{d})$ \reven{with} durations strictly greater than $\bar{d} - q$;
\item  $\mathcal{K}_2(\bar{d}, q) = \left\lbrace k \in \mathcal{K}(\bar{d}) ~ | ~ q \leq t_k \leq \bar{d} - q \right\rbrace$ \reven{contains} all the {\tour}s in $\mathcal{K}(\bar{d})$ \reven{with} durations greater than or equal to $q$ and lower than or equal to $\bar{d} - q$;
\item $\mathcal{K}_3(\bar{d},q) = \left\lbrace k \in \mathcal{K}(\bar{d}) ~ | ~  t_k < q \right\rbrace$ \reven{contains} all the {\tour}s in $\mathcal{K}(\bar{d})$ \reven{with} durations strictly lower than $q$.
\end{itemize}

\begin{proposition}
\label{prop:L2cuts}

The following constraints 
\begin{equation*}
\sum_{k \in \mathcal{K}_1(\bar{d},q)} \bar{d} \, \rho_{k} + \sum_{k \in \mathcal{K}_2(\bar{d},q)} t_k \, \rho_{k} \leq \bar{d} \times |\mathcal{P}| \;\;\;\; \forall
\bar{d} \in \mathcal{D}, q \in \mathcal{Q}(\bar{d})
 \qquad \constraintlabel{cnst:l2bound} \qquad \leftarrow \gamma^{\bar{d} \, q}_{MT} \leq 0 
\end{equation*}
are valid inequalities for $(M(\mathcal{K}))$.
\end{proposition}

\noindent\emph{Sketch of the proof.} For a detailed proof, see Appendix~\ref{proof:L2cuts}.
The proposition is proven by applying Proposition~\ref{prop:DFF} to Constraints~(\ref{cnst2BP}) with DFF $f_0^{\lambda}$, $\lambda \in \left[0;\frac{1}{2}\right]$:
\begin{align*}
f_0^{\lambda} : \left[ 0,1\right] &\rightarrow \left[0,1\right] \\
x & \mapsto \left\{
\begin{array}{l}
  1, \;\; \text{if } x > 1-\lambda \\
  x, \;\; \text{if } \lambda \leq x \leq 1-\lambda \\
  0, \;\; \text{if } x < \lambda
\end{array}
\right.
\end{align*}
\textcolor{red}{
$\hfill \square$\\
}





The meaning of these valid inequalities is the following.
Each picker $p \in \mathcal{P}$ is considered as a bin of size $\bar{d}$.
The {\tour}s in $\mathcal{K}_1(\bar{d}, q)$ of size strictly greater than $\bar{d} - q$ are accounted with size $\bar{d}$.
The {\tour}s of medium size, \ie those in $\mathcal{K}_2(\bar{d}, q)$ which lengths are between $q$ and $\bar{d} -q$, are accounted with their normal size $t_k$. 
Finally, the {\tour}s in $\mathcal{K}_3(\bar{d}, q)$ of size strictly lower than $q$ are ignored.

This is valid since a route of size greater than $\bar{d} - q$ and a route of size greater than $q$ cannot be performed by the same picker without exceeding the deadline $\bar{d}$.
Inequalities \eqref{cnst:l2bound} \reven{are} called Martello and Toth's cuts in the \reven{remainder} of the paper, and we denote by $\gamma^{\bar{d} \, q}_{MT} \leq 0$ the dual variable associated to the cut with parameters $\bar{d}$ and $q$.

\subsubsection{Fekete and Schepers's cuts}
\label{sec:FScuts}

Let us choose a deadline $\bar{d}\in\mathcal{D}$ and an integer $q\in\{2,\dots,\bar{d}\}$. For all integers $i\in\{0,\dots,q-1\}$ we define the following set:

$$
\mathcal{K}(\bar{d},i, q)\;=\;
\Bigl\{\;k\in\mathcal{K}: \bar{d}_k\leq\bar{d}, \text{ and } i\;\frac{\bar{d}}{q}<t_k\leq(i+1)\frac{\bar{d}}{q}\Bigr\}
$$

Consider the interval $[0,\bar{d}]$ divided in $q$ identical sub-intervals of length $\frac{\bar{d}}{q}$ and select a {\tour} $k\in\mathcal{K}(\bar{d},i,q)$. By definition of $\mathcal{K}(\bar{d},i, q)$ it holds $t_k-i\frac{\bar{d}}{q}>0$, \reven{therefore} route $k$ entirely covers $i$ intervals of $[0,\bar{d}]$, plus a strictly non-empty part of the $(i+1)$-th interval since, . 

It follows that given a deadline $\bar{d}\in\mathcal{D}$ and an integer $q\in\{2,\dots,\bar{d}\}$, we have at most $q-1$ equally sized intervals to place potential {\tour}s in $[0, \bar{d}]$, since the last interval is dedicated to allocate the residuals. Since this reasoning is valid for each picker $p\in\mathcal{P}$, it leads to the following proposition:

\vspace{1ex}
\begin{proposition}
\label{prop:FeketeCuts}

The following inequalities 
\begin{equation*}
\sum_{i=1}^{q-1} \sum_{k\in \mathcal{K}(\bar{d},i, q)} i \, \rho_k \leq (q-1) \times |\mathcal{P}|\quad\forall \bar{d}\in\mathcal{D}, q\in\{2,\dots,\bar{d}\}\quad \constraintlabel{eq:vi01} 
\qquad \leftarrow \gamma_{FS}^{\bar{d} \, q} \leq 0 
\end{equation*}
are valid for $(M(\mathcal{K}))$.
\end{proposition}

\noindent\emph{Sketch of the proof.} 
For a detailed proof see Appendix \ref{proof:FScuts}.
The proposition is proven by applying Proposition \ref{prop:DFF} to constraints~(\ref{cnst2BP}) with the DFF $g^{\lambda}(x)$, $\lambda \in \N \setminus \{0\}$:



\begin{align*}
g^{\lambda}(x) = & \left\{
\begin{array}{ll}
    0  & \text{if }  x = 0\\
  \frac{x(\lambda + 1) - 1}{\lambda} & \text{if }  x (\lambda +1) \in \Z\setminus\{0\} \\
  \left\lfloor (\lambda+1) x \right\rfloor \frac{1}{\lambda}, & \text{otherwise}
\end{array}
\right.
\end{align*}
\textcolor{red}{
$\hfill \square$\\
}

The DFF $g^{\lambda}(x)$ is a non-maximal DFF, meaning that a stronger cut may be obtained -- \reven{in particular with} the DFF proposed in \cite{Fekete2001}. However, the weaker version \reven{is easier to integrate into} the pricing problem  as will be detailed in Section~\ref{sec:vineqpricing}. 
Inequalities \eqref{eq:vi01} \reven{are} called Fekete and Schepers's cuts in the \reven{remainder} of the paper, and we denote by $\gamma_{FS}^{\bar{d} \, q} \leq 0$ the dual variable associated to the cut with parameters $\bar{d}$ and $q$.

\subsection{Valid inequalities for the JOBPRP}
\label{sec:cutsjobprp}

The problem described by formulation $(M(\mathcal{K}))$ is an extension of the JOBPRP \reven{with} packing constraints, \ie constraints~\eqref{cnst2:relax}.
Hence, the valid inequalities for the JOBPRP are also valid for the {\pbname}, thus valid for $(M(\mathcal{K}))$.

\subsubsection{Strengthened capacity cuts}
\label{sec:scccuts}
The strengthened capacity cuts (SCC) play with the size of the orders and the capacity of the trolley to impose the minimum number of {\tour}s needed to retrieve a given set of orders. They are defined as follows:

\begin{equation*}
    \sum_{k\in\mathcal K |\mathcal{O}_k \cap \mathcal{R} \neq \emptyset} \rho_k \geq \left\lceil \frac{\sum_{o \in \mathcal{R}}{\capa_o}}{\capa} \right\rceil\quad \forall \mathcal{R}\subseteq \mathcal O \qquad
\constraintlabel{eq:sr} \qquad  \leftarrow \gamma_{SC}^{\mathcal{R}} \geq 0
\end{equation*}

This family of cuts has been first proposed by \cite{Baldacci2008} for the Capacitated Vehicle Routing Problem (CVRP), and has also been used by \cite{BriantCCCLO20} and \cite{Wahlen2023} in the context of the JOBPRP.
Note that if there exists a set $\mathcal{R}' \subset \mathcal{R}$ such that $\left\lceil \frac{\sum_{o \in \mathcal{R}'}{\capa_o}}{\capa} \right\rceil = \left\lceil \frac{\sum_{o \in \mathcal{R}}{\capa_o}}{\capa} \right\rceil$, then the SCC defined over $\mathcal{R}'$ dominates the one defined over $\mathcal{R}$, and the latter does not need to be considered.
When $\mathcal{R} = \mathcal{O}$, the corresponding SCC states a constraint on the minimum number of routes required in any solution.

\subsubsection{Rank-1 cuts}
\label{sec:rank1}

The Chvátal–Gomory rank-1 cuts (R1C) are obtained by considering a small subset of the rows defining the set covering inequalities, \ie constraints~\eqref{cnst1:relax}.
Given a subset of orders $\mathcal{R} = \left\lbrace o_1, o_2, \ldots, o_q \right\rbrace \subseteq \mathcal{O}$, with $q \geq 3$ and non-negative multipliers $\textbf{p} = \left\lbrace p_1, p_2, \ldots, p_q \right\rbrace$, R1C are defined as:


\begin{equation*}
\sum_{k\in\mathcal K}\left\lfloor \sum_{o\in \mathcal{R} \cap \mathcal{O}_k} p_o \right\rfloor \rho_k \leq \left\lfloor \sum_{o \in \mathcal{R}}{p_o} \right\rfloor\quad \forall \mathcal{R} \subseteq \mathcal O, \textbf{p}\in\mathbb{R}_{+}^{q} \qquad
\constraintlabel{eq:rank1} \qquad  \leftarrow \gamma_{R1}^{\mathcal{R} \, \textbf{p}} \leq 0
\end{equation*}
%

The constraints obtained in the particular case when all multipliers have the same value $p_o = \frac{1}{h}$ for a given integer $h \in \left\lbrace 1, \ldots, |\mathcal{R}| \right\rbrace$ are known as subset-row inequalities.
First proposed by \cite{Jepsen2008}, \reven{they} are used in the state-of-the-art exact algorithms for vehicle routing problems \reven{because} they provide a significant improvement on the lower bound \citep{Pecin2017improved}.
R1Cs have been used by \cite{Muter2015} and \cite{Hessler2022} to solve the JOBPRP.

\cite{Pecin2017} performed a computational polyhedral study to determine the best possible vectors of multipliers $\textbf{p}$ for cuts with 3 to 5 rows.
For example, when $|\mathcal{R}| = 3$, there is a single optimal multiplier that is $\textbf{p} = \left\lbrace \frac{1}{2}, \frac{1}{2}, \frac{1}{2} \right\rbrace$.
The corresponding R1C is then:

\begin{equation*}
\sum_{k\in\mathcal K : |\mathcal{O}_k \cap \mathcal{R}| \geq 2}\rho_k \leq 1\quad \forall \mathcal{R} \subseteq \mathcal O, |\mathcal{R}| = 3 \qquad
\constraintlabel{eq:rank1f} 
\end{equation*}



%



%
%




\section{A column generation based heuristic to solve the \pbname}\label{sec:algo}

This section present\reven{s} the algorithm designed to tackle the {\pbname}. 
It is based on the one proposed in \cite{BriantCCCLO20} to solve the JOBPRP, and adapted to \reven{also handle} the sequencing aspect of the {\pbname}.
Briefly, it solves the linear relaxation $(LM(\mathcal{K}))$ of formulation $(M(\mathcal{K}))$ via column generation. $(LM(\mathcal{K}))$ is thus the so-called master problem, \reven{detailed} in Section~\ref{sec:master}.  Algorithm~\ref{alg:genColAlgo} \reven{gives} an overview of the complete algorithm.

\begin{algorithm}[h!t]
\caption{Column generation heuristic (\algo) pseudo code}
\label{alg:genColAlgo}
\begin{algorithmic}[1]
\State Computation of pool of {\tour}s $\pool$  (Section~\ref{sec:split})
\State $UB \leftarrow $ value of the best solution found by the initial heuristic (Section~\ref{sec:split})
\State Initialization of $\mathcal{K}'$ (Section~\ref{sec:initmaster})
\State $LB \leftarrow 0$
\Do
\While{$(LM)$ is not optimal and time limit is not reached} 
	\State Solve $(LM({\mathcal{K}}'))$ with columns in $\mathcal{K}'$ (Section~\ref{sec:master})
	\If{$\pool$ contains columns with negative reduced cost}\label{line:checkPool}
	    \State $\mathcal{K}' \leftarrow \mathcal{K}' \cup \{$negative reduced cost columns of $\pool\}$
	\Else
	    \State Solve pricing $Pr(\mathcal{K}^{'})$ (Section \ref{sec:pricing}) \label{lpr1}
	    \State $\mathcal{K}' \leftarrow \mathcal{K}' \cup \{$negative reduced cost columns found by $Pr(\mathcal{K}^{'})\}$
	    \State $LB \leftarrow \max\{LB, \textrm{ Lagrangian bound}\}$ (Section~\ref{sec:bound})
	\EndIf
	\State Compute $\mathcal{K}_{rich}$ from the most negative reduced cost column (Section~\ref{sec:pricingiter})\label{richsetl1}
    \State $\mathcal{K}' \leftarrow \mathcal{K}' \cup \mathcal{K}_{rich} $ \label{richsetl2}
	\State $\overline{UB} \leftarrow$ Seek for a primal solution from $\mathcal{K}_{rich}$
	\If{$\overline{UB} < UB$}
	    \State $UB \leftarrow \overline{UB}$
	 \EndIf
\EndWhile
\State Separate constraints \eqref{cnst:l2bound}, \eqref{eq:vi01}, \eqref{eq:sr} and \eqref{eq:rank1} (Section~\ref{sec:separation})
\doWhile{time limit not reached and at least one \eqref{cnst:l2bound}, \eqref{eq:vi01}, \eqref{eq:sr}, \eqref{eq:rank1} is violated}
\State $UB \leftarrow$ Solve $(B({\cal K}'))$
\end{algorithmic}
\end{algorithm}

\reven{The algorithm} first \reven{generates} a pool $\pool$ of {\tour}s and an initial solution  (Section~\ref{sec:split}). 
Then, \reven{it generates} an initial set of routes $\mathcal{K}'\subseteq\mathcal{K}$ over which the master problem is defined to obtain the restricted master problem $(LM(\mathcal{K}'))$ (Section~\ref{sec:initmaster}). 

The restricted master problem $(LM(\mathcal{K}'))$ is then solved. Based on the dual information, new negative reduced cost {\tour}s are included in $\mathcal{K}'$, first by scanning the pool $\pool$, and then by solving the pricing problem $(Pr(\mathcal{K}'))$ (Section~\ref{sec:pricing}) if no negative reduced cost {\tour} has been found in the pool. 
A rich column set $\mathcal{K}_{rich}$ is also included in $\mathcal{K}'$ from the most negative reduced cost column found either in the pool or by the pricing problem (Section~\ref{sec:pricingiter}).

After each iteration of the pricing problem, \reven{we compute} the Lagrangian bound \reven{to possibly update} the value of the lower bound $LB$ (Section~\ref{sec:bound}). Moreover, in the hope of improving the upper bound, some columns of the rich column set are generated to provide a feasible solution for $(B(\mathcal{K}))$.

To strengthen the lower bound provided by the resolution of formulation $(LM(\mathcal{K}))$, we consider the valid inequalities presented in the previous section: the Martello and Toth's cuts \eqref{cnst:l2bound}, the Fekete and Schepers's cuts \eqref{eq:vi01}, the strengthened capacity cuts \eqref{eq:sr} and the rank-1 cuts \eqref{eq:rank1}. Their separation is detailed in Section~\ref{sec:separation}.

The procedure continues until $(LM(\mathcal{K}))$ is optimally solved or the time limit is reached. A final resolution of formulation $(\mathcal{B}(\mathcal{K}'))$  defined over the set of {\tour}s $\mathcal{K}'$ generated so far is launched with a time limit $\tau_{l}$. 
This column generation based heuristic procedure is called {\algo} in the following.

Since the proposed algorithm is based on the one proposed in \cite{BriantCCCLO20}, \reven{we do not report here} the details about the solving of the pricing problem, the rich column set and the Lagrangian bound. \reven{Instead, we}
concentrate on and highlight the elements that allow to extend the algorithm of \cite{BriantCCCLO20} to solve the {\pbname}.
In particular, the main contributions are:
\begin{itemize}
    \item a Mixed Integer Program (MIP) formulation of the pricing problem to take into account the bin packing constraints~\eqref{cnst2:relax} (Section~\ref{sec:mippricing});
    \item the strengthening of the tour constraints used in the MIP formulation of the pricing problem (Section~\ref{sec:tourCst});
    \item the impact of the Martello and Toth's cuts \eqref{cnst:l2bound}, and Fekete and Schepers's cuts \eqref{eq:vi01} on the MIP formulation of the pricing problem (Section~\ref{sec:vineqpricing}), and their separation (Section~\ref{sec:separation});
    \item the determination of an initial solution for the {\pbname} and a pool of columns (Section~\ref{sec:split}).
\end{itemize}
Contrary to \cite{BriantCCCLO20} the proposed algorithm \reven{uses no stabilization}.
Indeed, preliminary experimental results detailed in Appendix~\ref{annex_stabilization} show that adding stabilization as proposed in \cite{BriantCCCLO20} does not improve the quality of the results.

The next sections present in detail each component of the procedure. 


\subsection{Master problem}\label{sec:master}

The master problem that we solve across column generation is the linear relaxation of formulation $(M(\mathcal{K}))$, denoted as $(LM(\mathcal{K}))$, where the integrality requirements on variables $\rho_{k}, k\in\mathcal{K}$ are removed. The master problem $(LM(\mathcal{K}))$ reads as follows:
$$
(LM(\mathcal{K}))\left\{
\begin{array}{rclllllllllllll}
\min&\displaystyle\sum_{k\in\mathcal{K}}t_k\;\rho_{k}&&&&\constraintlabel{obj:master}\\
\\
s.t.
&\displaystyle\sum_{k\in\mathcal{K}(o)}\;\rho_{k}\;\geq\; 1
&&\forall o\in\mathcal{O}&&\constraintlabel{cnst1:master}&\leftarrow\;\alpha^o\geq 0\\
\\
&\displaystyle\sum_{k\in\mathcal{K}(\bar{d})}t_{k}\;\rho_{k}\;\leq\; \bar{d}\times |\mathcal{P}|
&&\forall \bar{d}\in\mathcal{D}&&\constraintlabel{cnst2:master}&\leftarrow\;\beta^{\bar{d}}\leq 0\\
\\
&\rho_{k}\;\geq 0
&&\forall k\in\mathcal{K}&&\constraintlabel{vardef:master}\\
\\
\end{array}
\right.
$$

The dual variables associated with each family of constraint\reven{s} are reported next to the expression. Usually, the master problem is solved over a restricted set of {\tour}s $\mathcal{K}'\subseteq\mathcal{K}$. We thus refer to $(LM(\mathcal{K}'))$ as the restricted master problem, as is usually done in the literature. 

\subsection{Pricing problem}\label{sec:pricing}

In practice, the master problem $(LM(\mathcal{K}))$ is solved on a subset $\mathcal{K}^{'} \subseteq \mathcal{K}$ of the entire set of routes. The role of the pricing problem is to identify a new variable among those in $\mathcal{K}\setminus\mathcal{K}^{'}$ of negative reduced cost, or to prove that none exists. The reduced cost $\bar{c}_k$ of variable $\rho_{k} \in \mathcal{K}$ is computed as follows: 
$$
\bar{c}_k\;=\;t_k\;
-\;\sum_{o\in\mathcal{O}}\alpha^o\;
\mathbb{1}_{\mathcal{K}(o)}(k)\;
-\;\sum_{\bar{d}\in\mathcal{D}}\beta^{\bar{d}}\;t_k\;
\mathbb{1}_{\mathcal{K}(\bar{d})}(k)
$$

In the rest of the section, we first formulate the pricing problem as a MIP presented in Section~\ref{sec:mippricing}. Sections~\ref{sec:tourCst} and~\ref{sec:coeffs} present the strengthening of the tour constraints of~\cite{BriantCCCLO20} that are dynamically included in the MIP of the pricing.
Section~\ref{sec:vineqpricing} discusses the impact of the valid inequalities on the modelling and resolution of the pricing problem.

\subsubsection{A MIP formulation for the solving of the pricing problem}\label{sec:mippricing}

The solving of the pricing problem is performed as  in~\cite{BriantCCCLO20}, \ie \reven{with} a cutting plane algorithm.
In this section we provide details on the MIP formulation of the pricing problem; details about the solving of the pricing can be found in Section~\ref{sec:pricingiter}.
For \reven{the} sake of simplicity, we state the MIP formulation without considering the dual contribution deriving from  valid inequalities~\eqref{cnst:l2bound},~\eqref{eq:vi01},~\eqref{eq:sr} and~\eqref{eq:rank1} presented in Section~\ref{sec:ineqdff}. This topic is covered \reven{separately} in Section~\ref{sec:vineqpricing}. 

The MIP formulation of the pricing problem presented here extends the one in~\cite{BriantCCCLO20} by taking into account constraints~\eqref{cnst2:master} of the master problem which capture the additional sequencing characteristic of the {\pbname}. 
The decision variables are:

\begin{itemize}
\item $t\geq 0$: real non-negative variable that represents a lower bound on the total time of the route (including picking and setup times);
\item $\delta^{\bar{d}} \geq 0$: real non-negative variable equal to $t$ if the deadline of the route is less than $\bar{d}$, 0 otherwise;
\item $e^o \in \{0,1\}$: binary variable equal to 1 if order $o$ is in the route, 0 otherwise;
\item $\mu^{\bar{d}} \in \{0,1\}$: binary variable equal to 1 if the deadline of the route equals $\bar{d}$, 0 otherwise.
\end{itemize}

We \reven{denote as} $\mathcal{O}^{\bar{d}} = \left\lbrace o \in \mathcal{O} : \bar{d}_o = \bar{d} \right\rbrace$ the set of orders with deadline $\bar{d}$. 

The pricing is solved by a cutting plane algorithm based on the following MIP model $Pr(\mathcal{K}^{'})$:

$$
Pr(\mathcal{K}^{'})\;\left\{
\begin{array}{rcllllllllll}
\min
&\displaystyle
\;t\;-\;\sum_{o\in\mathcal{O}}\alpha^o\;e^o\;-\;\sum_{\bar{d}\in\mathcal{D}}\beta^{\bar{d}}\;\delta^{\bar{d}}
&&\constraintlabel{pr:obj}\\\\
s.t.
&\displaystyle
\sum_{o\in\mathcal O}\capa_o e_o \leq \capa
&&\constraintlabel{pr:01}
\\\\
&\displaystyle
t_k \left( \sum_{o\in\mathcal{O}_k} e_o - |k| + 1 \right) \leq t
&\forall k\in\mathcal{K}' 
&\constraintlabel{cst:cutpricing}
\\\\
&\displaystyle\sum_{\bar{d}\in\mathcal{D}}\mu^{\bar{d}}\;=\; 1
&&\constraintlabel{pr:03}\\
\\
&\displaystyle\sum_{o\in\mathcal{O}^{\bar{d}}}e^o\;\geq\; \mu^{\bar{d}}
&\forall \bar{d}\in\mathcal{D}&\constraintlabel{pr:04}\\
\\
&\displaystyle\sum_{\bar{d}\leq\bar{d}_o}\mu^{\bar{d}}\;\geq\;e^o
&\forall o\in\mathcal{O}&\constraintlabel{pr:05}\\
\\
&\displaystyle\delta^{\bar{d}}\;\geq\;t\;-\;\sum_{\bar{d}'> \bar{d}}\bar{d}'\mu^{\bar{d}'}
&\forall \bar{d}\in\mathcal{D}&\constraintlabel{pr:06}\\
\\
&\displaystyle t\;\leq\;\sum_{\bar{d}\in\mathcal{D}}\bar{d}\;\mu^{\bar{d}}
&&\constraintlabel{pr:07}\\
\\
&e^o\;\in\;\{0,1\}
&\forall o\in\mathcal{O}
&\constraintlabel{pr:08}\\
&\mu^{\bar{d}}\;\in\;\{0,1\}
&\forall \bar{d}\in\mathcal{D}
&\constraintlabel{pr:09}\\
&\delta^{\bar{d}}\;\geq\;0
&\forall \bar{d}\in\mathcal{D}
&\constraintlabel{pr:10}\\
&t\;\geq\;0&&\constraintlabel{pr:11}\\
\end{array}
\right.
$$

The objective function \eqref{pr:obj} of $Pr(\mathcal{K}^{'})$ minimizes the reduced cost expression.
Constraint \eqref{pr:01} enforce\reven{s} the respect of the capacity of the trolley.
Variable $t$ encodes a lower bound on the travel time of a route. It is strengthened by the dynamic generation of constraints \eqref{cst:cutpricing} that are called the {\em tour constraints}. Typically, any selection of orders that is a superset of a known set of orders $k$ require\reven{s} a total time longer than $t_k$\rev{,} so that $t \geq t_k$.
Constraint \eqref{pr:03} state\reven{s} that a single deadline \reven{is} selected for the route. 
Constraints \eqref{pr:04} impose that at least one order in the route \reven{has} the selected deadline. 
Constraints \eqref{pr:05} impose that the deadline of the route is smaller than the \reven{deadline} of any order picked by the route. Note that \reven{together, those} two previous constraints \eqref{pr:04} and \eqref{pr:05} ensure that the deadline of the route is the minimum among the deadlines of all picked orders.
Constraints \eqref{pr:06} guarantee that the $\delta^{\bar{d}}$ variables are consistent with $\mu^{\bar{d}}$. 
Constraint \eqref{pr:07} impose\reven{s} that the duration of the route is smaller than its deadline to avoid obviously inconsistent routes.
Finally, constraints \eqref{pr:08}--\eqref{pr:11} define the domain of the variables.

It is important to notice that solving $Pr(\mathcal{K}')$ provide\reven{s} a lower bound on the value of the lowest reduced cost route in $\mathcal{K}$, since constraints~\eqref{cst:cutpricing} are only defined for routes in set $\mathcal{K}'$.
When solving $Pr(\mathcal{K}')$, the value of variable $t$ also represent\reven{s} a lower bound on the exact time needed to collect all the orders selected by the pricing problem ($e_o=1$).

\subsubsection{Tour constraints and their strengthening}\label{sec:tourCst}

As mentioned in the previous section, to enhance the resolution of the pricing problem, \cite{BriantCCCLO20} propose to add to $Pr(\mathcal{K}')$, for each route $k \in \mathcal{K}'$, the following so-called {\em tour constraint}:
$$
    t_k  \left( \sum_{o \in\mathcal{O}_k}{e_o} - |k| + 1 \right) \leq t \qquad \forall k \in \mathcal{K}' \qquad \eqref{cst:cutpricing}
$$

The rationale behind these constraints is the following: any selection of orders that is a superset of a known set of orders $k$ require\reven{s} a total time longer than $t_k$, so that $t \geq t_k$.
Moreover, a tour constraint is active only when all variables $e_o, o \in\mathcal{O}_ k$ involved in it take value 1. Otherwise, the term between brackets is lower than or equal to zero and the constraint is inactive.

We propose to strengthen these constraints by \reven{focusing on} the times that account for the total time $t_k$ of a route. These are the setup time $t^{setup}$, the picking times $t_o^{pick}$ for each order that is picked, and the travelling times.

Let us define generic coefficients $b\geq0$ and $a_o\geq0, \forall o \in \mathcal{O}$.
Let us also assume that the following constraint is valid:
\begin{equation*}
    b  + \sum_{o \in\mathcal{O}_k} a_o \,e_o \leq t
    \qquad \constraintlabel{eq:eq2}
\end{equation*}

Let us \reven{de}note 
$$\displaystyle\Delta_k=t_k - b - \sum_{o \in\mathcal{O}_k} a_o$$
as the difference between $t_k$ and the left-hand side of constraint~\eqref{eq:eq2} evaluated on {\tour} $k$, \ie when $e_o=1$ for all $o \in\mathcal{O}_k$. 
Since $t$ encodes in $Pr(\mathcal{K}')$ a lower bound on the travel time, it holds $t \leq t_k$ where $k$ is the optimal route found by solving $Pr(\mathcal{K}')$. 
Note that $t = t_k$ is insured if $k \in \mathcal{K}'$. Thus, \reven{assuming} constraints~\eqref{eq:eq2} \reven{is valid}, it follows that $\Delta_k \geq 0$.

The following propositions hold.

\begin{proposition}
If constraints~\eqref{eq:eq2} are valid, then the following constraints 
\begin{equation*}
    b  + \sum_{o \in\mathcal{O}_k} a_o \,e_o \;+\; \Delta_k \left( \sum_{o \in\mathcal{O}_k} e_o - |k| + 1 \right) \leq t
    \qquad \forall k \in \mathcal{K}' \qquad \constraintlabel{cst:cutpricing2}
\end{equation*}
are valid for $Pr(\mathcal{K}')$ and allow to exactly compute the value of $t_k$ when $e_o=1$ for all orders $o \in \mathcal{O}_k$.
\end{proposition}

\begin{proof}
If $e_o=1$ for all $o\in\mathcal{O}_k$, then constraints~\eqref{cst:cutpricing2} applied to route $k$ become $t_k\leq t$.
Otherwise, we have: $\displaystyle\sum_{o \in\mathcal{O}_k} e_o - |k|+1 \leq 0$, and the term added to constraint~\eqref{eq:eq2} to obtain constraints~\eqref{cst:cutpricing2} is non-positive since $\Delta_k\geq 0$.
\end{proof}


\begin{proposition}
Constraints \eqref{cst:cutpricing2} provide a better estimation than constraints \eqref{cst:cutpricing} of the total time of the route determined by the resolution of the pricing problem $Pr(\mathcal{K}')$.
\end{proposition}

\begin{proof}
Let us suppose that the solution of $Pr(\mathcal{K}')$ provides a {\tour} $k^*$, and let us consider a route $k\in\mathcal{K}'$.
If $k^* \cap k = k$, \ie $k^*$ is a superset of $k$, then constraints~\eqref{cst:cutpricing} and~\eqref{cst:cutpricing2} both lead to $t_k \leq t$.

If $k^* \cap k \subsetneq k$, then the following relation holds:
$$
t 
\;\;\;\;\geq\;\;\;\;
\underbrace{
b  + \sum_{o \in\mathcal{O}_{k^*}} a_o \;+\; \Delta_k (|k^* \cap k| - |k| + 1)
}_{\mbox{Constraint \eqref{cst:cutpricing2}} }
\;\;\;\;\geq\;\;\;\;
\underbrace{
t_k  (|k^* \cap k| - |k| + 1).
}_{\mbox{Constraint \eqref{cst:cutpricing}}}
$$
Indeed, $\displaystyle b+\sum_{o\in\mathcal{O}_{k^*}}a_o\geq 0$, $\Delta_k\leq t_k$ and $(|k^* \cap k| - |k| + 1)\leq 0$. 
Then $t$ takes a value that is closer to $t_{k^*}$ when using constraints~\eqref{cst:cutpricing2} instead of constraints~\eqref{cst:cutpricing}.
\end{proof}

It is possible to generalize constraints~\eqref{cst:cutpricing2} to consider all orders $o \in \mathcal{O}$. 
To this end, \reven{we} consider a route $k$ and assume that the following constraint \reven{is valid}:
\begin{equation*}
    b  + \sum_{o \in \mathcal{O}} a_o \,e_o \leq t
    \qquad \constraintlabel{eq:eq2b}
\end{equation*}
Lt us consider coefficients $a_o$ such that:
\begin{equation*}
    t_k + \sum_{o\in \mathcal{O}_{k'}\setminus\mathcal{O}_k} a_o \leq t_{k'}\;\;\;\;\forall k'\supseteq k 
    \qquad \constraintlabel{eq:eq2bhyp}
\end{equation*}

We thus state the following proposition.
\begin{proposition}
\label{prop:validPricingCuts}
If constraint~\eqref{eq:eq2b} is valid and hypothesis~\eqref{eq:eq2bhyp} holds, then the following constraints
\begin{equation*}
    b  + \sum_{o \in \mathcal{O}} a_o \,e_o \;+\; \Delta_k \left( \sum_{o \in\mathcal{O}_k} e_o - |k| + 1 \right) \leq t
    \qquad \forall k \in \mathcal{K}' \qquad \constraintlabel{cst:cutpricing3}
\end{equation*}
are valid for $Pr(\mathcal{K}')$ and allow to compute $t_k$ exactly.
\end{proposition}

\begin{proof}
See Appendix \ref{proof:validPricingCuts}.
\end{proof}

Note that, since constraints~\eqref{cst:cutpricing2} \reven{can easily be seen as} weaker than constraints~\eqref{cst:cutpricing3}, we only add the latter to the pricing problem.



\subsubsection{Examples of coefficients to be used in constraints~\eqref{cst:cutpricing3}}\label{sec:coeffs}
By setting $b=t^{setup}$, and $a_o=t^{pick}_o$ for all $o \in \mathcal{O}$, it is clear that constraint~\eqref{eq:eq2b} is satisfied, and hypothesis~\eqref{eq:eq2bhyp} holds.
Such coefficients can be easily used in constraint~\eqref{cst:cutpricing3}.
\reven{They reflect} that the duration $t_k$ of a route $k$ is not only related to travelling time, but also to setup and picking times.

\reven{We can then} provide larger coefficients $a_o$ for orders $o \in k$, by considering the minimum additional traveling time to pick each order $o \in k$.
To this end, given a sequence $(o_1, \ldots, o_{|k|})$ of the orders in route $k$, we define:
\begin{equation*}
a_{o_i} = t^{pick}_{o_i} + \min_{\mathcal{R} \subseteq \{ o_1, \ldots, o_{i-1} \}}{\left\lbrace \tilde{t}_{\mathcal{R} \cup o_i} - \tilde{t}_{\mathcal{R}} \right\rbrace} \qquad \forall i \in (1, \ldots, |k|)
\end{equation*}
where $\tilde{t}_{\mathcal{R}}$ is the optimal traveling time to collect all orders in $\mathcal{R}$, and $\tilde{t}_{\emptyset} = 0$.
Note that we still define $b = t^{setup}$ and $a_o=t^{pick}_o$ for all $o \in \mathcal{O} \setminus k$.

Coefficient $a_{o_i}$ for each order $o_i \in k$ corresponds to the picking time plus the minimum additional travelling time required when adding order $o_i$ to any subset of orders in $\{ o_1, \ldots, o_{i-1} \}$.
However, note that the computation of each coefficient has a complexity of $O(2^i)$.

In order to avoid large computation times to compute the coefficients when the constraint considers a large number of orders, we propose to compute a lower bound on these coefficients.
The latter is based on the computation of an estimation of the travelling time to collect a set of orders and is calculated by considering the rectangular bounding box of the locations to be visited.
Details are provided thereafter.

Given a subset of orders $\mathcal{R} \subseteq \mathcal{O}$, \reven{we} denote by $\tilde{t}^{box}_{\mathcal{R}}$ the time required to travel the perimeter of the 2-dimensions bounding box of the locations of all items in $\mathcal{R}$ and the depot.
\reven{Interestingly,} an optimal route can go outside the bounding box if the item in $\mathcal{R}$ with the largest $y$ coordinate is located in the second half of its block (\ie nearest from the upper cross aisle of the block).
In such a case, an optimal route could use the upper cross-aisle of the upper block.
We name {\em large bounding box} the bounding box of the locations of all items in $\mathcal{R}$ and the depot, enlarged to the upper cross-aisle.
We denote by $\tilde{t}^{boxLarge}_{\mathcal{R}}$ the time required to travel the perimeter of the large bounding box.
Then, given a sequence $(o_1, \ldots, o_{|k|})$ of the orders in route $k$, the following coefficients can be used:
\begin{equation*}
a_{o_i} = t^{pick}_{o_i} + \max{ \left\lbrace 0 ; \tilde{t}^{box}_{\{ o_1, \ldots, o_{i} \}} - \tilde{t}^{boxLarge}_{\{ o_1, \ldots, o_{i-1} \}} \right\rbrace }
\end{equation*}

Coefficient $a_{o_i}$ is equal to the picking time plus a lower bound on the additional travelling time required when adding order $o_i$ to any subset of orders in $\{ o_1, \ldots, o_{i-1} \}$.

\subsubsection{Impact of valid inequalities on the pricing problem}\label{sec:vineqpricing}

Column generation-based algorithms are usually used to solve a problem when its mathematical formulation is defined over a set 
with an exponential amount of variables.
For real-size instances, the resulting number of variables does not allow 
to generate them all a-priori, 
thus only the interesting ones are dynamically generated.
\reven{The so-called \emph{extended} formulations are} obtained by applying Dantzig-Wolfe reformulation techniques to formulations for the same problem defined over a polynomial number of variables. Such formulations are usually called {\em compact}.

Cutting planes added to the extended formulation are called {\em robust} cuts when they do not increase the complexity of the pricing problem needed to generate interesting variables on the fly. This is usually the case of constraints that can be expressed with the variables of the compact formulation. Opposite to that, when the complexity of the pricing problem  increase\reven{s} by \reven{considering} a family of cutting planes, they are {\em nonrobust}. Hence, each nonrobust cut added to the extended formulation (\ie the master problem) directly makes the pricing problem more difficult to solve. On the opposite, nonrobust cuts are known for their great potential for reducing integrality gaps \citep{contardo2019}.

The four families of cuts proposed in Sections~\ref{sec:ineqdff} and~\ref{sec:cutsjobprp} are \emph{nonrobust} cuts. 
When the pricing problem is modelled as an elementary shortest path problem with resource constraints, as usually done for routing problems, and solved by dynamic programming via a labelling algorithm, the management of non-robust cuts can be cumbersome \citep{Pecin2017improved}.
However, when the pricing problem is modelled as a MIP, as \reven{in} the present work, and solved with a commercial solver, adding nonrobust cuts in the master problem implies adding extra decision variables and constraints in the MIP formulation. \reven{Although it can enlarge} the MIP and increase its resolution time, the coding effort to manage them is rather limited.
Each family of cuts therefore requires to update the MIP formulation of the pricing problem. \reven{These} modifications are detailed in Section \ref{annex_pricing} of the Appendix since \reven{they} amount to relatively standard modelling techniques in integer programming.

\subsection{Pricing iteration}
\label{sec:pricingiter}

After the resolution of the restricted master problem $(LM(\mathcal{K}'))$, we first check if 
pool $\pool$ (see Section~\ref{sec:split}) includes {\tour}s \reven{of} negative reduced cost. 
If any, they are added to $\mathcal{K}'$. 
If not, the pricing problem is solved. 
As explained in Section~\ref{sec:pricing}, solving $Pr(\mathcal{K}')$ provides a lower bound on the most negative reduced cost column in set $\mathcal{K}$.
When solving $Pr(\mathcal{K}')$, each time a feasible route with a negative value is found, such a route can potentially lead to a negative reduced cost column.
An efficient dynamic programming algorithm (see \citet{Cambazard18, Pansart18}) is then applied on such promising routes to compute their exact total time and check whether they indeed correspond to negative reduced cost routes.

After solving $Pr(\mathcal{K}')$, it is possible that no real negative reduced cost route has been found, while the optimal value of $Pr(\mathcal{K}')$ was negative.
In such a case, set $\mathcal{K}'$ is enlarged with the best route found by $Pr(\mathcal{K}')$, and $Pr(\mathcal{K}')$ is solved again.
More details on this iterative solving of $Pr(\mathcal{K}')$ can be found in \citet{BriantCCCLO20}.
Moreover, in this work, $Pr(\mathcal{K}')$ is solved with a time limit $\tau_{Pr}$.
When solving the pricing problem \reven{yields} at least one negative reduced cost route, we add to $\mathcal{K}'$ up to 10 negative reduced cost {\tour}s found during the resolution of $Pr(\mathcal{K}')$.

After finding negative reduced cost routes, either in the pool or by solving $Pr(\mathcal{K}')$, we then determine a set of {\tour}s to add in $\mathcal{K}'$ with respect to the {\tour} with the most negative reduced cost. 
This route will be denoted as $k^*$ in the following.
As proposed in \cite{BriantCCCLO20}, two different strategies are considered to \reven{determine} other {\tour}s to insert in $\mathcal{K}'$: the first strategy is to complete column $k^*$ with a set of columns that potentially constitute a feasible solution of $(B(\mathcal K))$, and the second strategy is to \reven{use} columns similar to $k^*$ with fewer orders.
For the first strategy, we proceed as in \cite{BriantCCCLO20} to generate a set of promising columns such that each order is collected.
However, because of the bin packing constraints~\eqref{cnst2BP}, there is no guarantee on the feasibility of the solution.
Hence, we then solve formulation $(B(\mathcal{K}))$ with this set of promising columns to evaluate the feasibility of the solution.
If the solution is feasible and improves the current upper bound, then  the upper bound is updated accordingly.
For the second strategy, if $k^*$ has less than 8 orders, as in \cite{BriantCCCLO20}, we generate all sub-tours of $k^*$ that use no more than 75\% of the capacity of the trolley.
If $k^*$ has 8 orders or more, we generate all sub-tours of $k^*$ with one order less than $k^*$.
All columns generated by the two strategies compose the rich column set $\mathcal{K}_{rich}$ mentioned in  Algorithm~\ref{alg:genColAlgo}. 
The reader is referred to \cite{BriantCCCLO20} for \reven{further} details.

\subsection{Separation of the valid inequalities}\label{sec:separation}

Martello and Toth's cuts \eqref{cnst:l2bound} can be stronger than the relaxation of the bin packing constraints~\eqref{cnst2:relax} only if $\mathcal{K}_1(\bar{d}, q)$ contains at least one element $k$ such that $\rho^*_k > 0$, where $\rho^*$ is the optimal solution of the master problem $(LM(\mathcal{K}))$.
Note that when $q$ increases, the size of $\mathcal{K}_1(\bar{d},q)$ increases, but the size of $\mathcal{K}_1(\bar{d},q) \cup \mathcal{K}_2(\bar{d},q)$ decreases.
So, \reven{no} value of $q$ \reven{can} a priori provide the most violated Martello and Toth's cut~\eqref{cnst:l2bound}.
However, as mentioned in \cite{Martello1990}, it is not necessary to consider all possible values of $q$ in $\mathcal{Q}(\bar{d})$.
The size of $\mathcal{K}_1(\bar{d},q)$ increases when there exist\reven{s} $k \in \mathcal{K}(\bar{d})$ such that $q = \bar{d} - t_k +1$.
For a given size of $\mathcal{K}_1(\bar{d},q)$, the most violated Martello and Toth's cut \eqref{cnst:l2bound} is the one with the lowest $q$ value, \ie associated with the highest size of $\mathcal{K}_2(\bar{d},q)$.
\reven{Following} all these remarks, for a given $\bar{d}$, we just consider the $q = \bar{d} - t_k +1$ values such that there is a route $k$ in $\mathcal{K}(\bar{d})$ with $\rho^*_k > 0$ and $t_k \geq \frac{1}{2}\bar{d} + 1$.

For the Fekete and Schepers's cuts \eqref{eq:vi01}, given a value of $\bar{d}$, several inequalities can be generated with different values of $q$.
We decided to \reven{include in the model} the cut with the highest violation, \reven{following} a preliminary computational study \reven{that pointed to this configuration as the best}.

For the SCC cuts \eqref{eq:sr}, as done in \cite{BriantCCCLO20}, we only consider those characterized by a right-hand side that equals the minimum number of trolleys to retrieve all orders. 
Note that given the benchmark of instances, it is possible in practice to enumerate all subsets of orders that require this minimum number of trolley to be picked. As a consequence, we do not call a separation algorithm, but scan the set of such subsets to look for violated constraints. 
Moreover, we only consider the SCC cuts that are minimum for inclusion.
The master problem is initialized only with the SCC cuts \eqref{eq:sr} defined over the entire set of orders. The others are added to the master problem only if they are violated. 

For the rank-1 cuts \eqref{eq:rank1}, after each resolution of the master problem, we check, in this order, if such constraints defined over subset\reven{s} of order size 3, then size 4 are violated. 
For small size instances (with 18 orders or less), we also check subset\reven{s} of order size 5.
We stop the procedure when 40 violated constraints are found or when the list has been fully scanned. In any case, we do include in the master problem at most the 20 most violated cuts.

\subsection{Computation of a lower bound}\label{sec:bound}

When the procedure reaches the time limit before $(LM(\mathcal{K}))$ is solved, we calculate the so-called Lagrangian bound to be able to provide a valid lower bound on the value of the optimal solution. The interested reader is referred to \cite{BriantCCCLO20} for a detailed explanation of its computation. 

\subsection{Initial solution and generation of a pool of {\tour}s}\label{sec:split}

Before starting the resolution of the master problem via column generation, we provide an initial solution and generate a set of possible promising {\tour}s to include in a pool $\pool$. 

We first create an initial solution using a greedy procedure. The first step consists in creating one batch for each order. Then, the procedure sorts batches by increasing deadlines and then assigns batches to pickers in the following way. Pickers are arbitrarily ranked and the first batch is assigned to the first picker, the second batch to the second picker, and so on until all pickers have been assigned \reven{to} a batch. Then the next batch is assigned again to the first picker, the second next batch to the second picker and so on until all batches have been assigned.

It is likely that the initial solution \reven{obtained through this greedy procedure} violates some deadlines. For this reason\rev{,} we allow \rev{unfeasibility}, which is measured as the tardiness of the solution. The tardiness is computed as $T = \sum_{p \in \mathcal{P}} T_p$ where $T_p = \sum_{r_k \in \mathcal{K}_p} \max(0, t_k - \overline{d}_k)$ is the tardiness of the picker computed from his/her batches or routes. The tardiness of a route is computed from a single order, the one with the minimum deadline (recall that $\overline{d}_k = \min_{o \in \mathcal{O}_k} \overline{d}_o$).

\reven{With} the aim of recovering feasibility or improving the quality of the solution, \reven{this initial solution} then undergoes an improvement phase \reven{with} a local search (LS) procedure composed \reven{of} the following operators: 
\begin{itemize}
    \item relocation of an order in a batch of another picker;
    \item relocation of an order in a different batch of the same picker;
    \item exchange of orders from different batches of the same picker;
    \item exchange of orders from different batches of different pickers;
    \item relocation of a batch from one picker to another;
    \item exchange of batches between two different pickers. 
\end{itemize}
The sequence of batches/routes $r_k$ for each picker is maintained sorted by increasing deadline $d_k$ (when two routes have the same deadline, the shortest one is put first). \reven{This comes from} the bin packing view of the problem: the local search does not operate on the sequence itself by moving/swapping batches of the same picker.
\reven{Sorting} the routes by deadlines \reven{is enough} to determine the feasibility of the solution, \ie \reven{whether it incurs} tardiness or not.

At each iteration of the LS the neighborhoods defined by each operator are fully explored. Once the exploration of one neighbourhood is \reven{completed}, the best move, if it improves the solution, is applied. The next neighbourhood is then explored. The neighbourhoods are explored in the order they appear in the previous list. This order is arbitrary.

\reven{Once} the neighbours have been explored, the route duration of all batches in the solution is optimized by means of the dynamic programming procedure described in \cite{Cambazard18}. This provides the exact time to collect all the orders in a batch. If at least one operator has improved the solution, the procedure is run again, otherwise it is stopped.

We consider that a move improves a solution if one of the following two conditions apply:
\begin{itemize}
    \item the solution \reven{incurs tardiness that is reduced by the move;} 
    \item the move does not reduce nor increase the tardiness of the solution, but reduces the total time needed to collect all the orders.
\end{itemize}

To increase the chance \reven{for} the procedure to obtain a high-quality solution, the best solution obtained so far is perturbed by randomly removing and reinserting all orders from a number of batches. Typically, a batch $b$ is randomly selected and each order of $b$ is moved to another existing batch if the capacity allows it. At most three target batches are tried for each order and orders of $b$ that are left (\ie that have not been successfully inserted somewhere else) form a batch assigned to a randomly selected picker. Note that this procedure may fail to perturb a given solution if the batches are all \reven{at} full capacity (or in general if the residual capacity of batches does not allow any order to be inserted). \reven{In} this situation, the selected batch is replaced by as many batches as the orders it contains. 
The new solution undergoes the LS. 
If a new best solution is obtained, the number of batches to remove in the perturbation phase is set to a minimum value, otherwise it is increased by one. 
The whole procedure stops after a certain time limit $\tau_{pool}$ is reached.

In order to obtain an initial pool $\pool$ of {\tour}s to feed the master problem $(LM(\mathcal{K}'))$, all solutions obtained after a call to the LS are stored in memory. 
To avoid keeping duplicates we discard solutions with the same tardiness and cost \reven{than} another already present in the memory. 
\reven{The pool of routes is then built including all the {\tour}s of each solution in the memory.}

Note that the heuristic used to generate a pool of {\tour}s is designed to quickly provide a set of diverse {\tour}s.
The local search component is thus dedicated to obtain a \rev{good-quality} solution fast, based on a newly perturbed solution -- not to optimize a \rev{good-quality} solution.

\subsection{Initialization of the restricted master problem}\label{sec:initmaster}

The restricted master problem is initialized by adding in $\mathcal{K}'$ all {\tour}s of size 1. Then $\mathcal{K}'$ is completed with at most 2000 {\tour}s of size 2, 3 and 4 that fill the trolley by at least 90\% of its capacity, starting by routes of size 2, then size 3 and finally size 4. 
Then, for all the \reven{generated} tours \reven{that} do not use all the capacity, we add some orders if possible to fill the remaining capacity.
We also add in $\mathcal{K}'$ the {\tour}s \reven{from} the best solution found by the procedure described in Section~\ref{sec:split}.

\reven{Regarding} the valid inequalities, from the beginning we only add the SCC cut defined over the entire set of orders $\mathcal O$.

\section{Experimental results}
\label{se:expes}
The experiments were performed on an Intel(R) Core(TM) i7-8650U CPU @ 2.11 GHz processor with 16~GB of RAM and each algorithm ran 
on a \rev{single-thread}. The code was written in Java and CPLEX 20.1 was used to solve the LPs and the MIPs.
The time limit to solve the small instances is set to 3 600 seconds while the time limit to solve the large instances is set to 14 400 seconds (4 hours).



\subsection{Benchmark of instances}

To evaluate the performance of procedure {\algo} proposed in Section~\ref{sec:algo}, we use the benchmark of instances generated to evaluate the iterated local search (ILS) algorithm proposed in~\cite{VanGils2019}. This benchmark of instances is divided \rev{into} two sets: the \emph{small instances} made of 6, 12 and 18 orders, and the \emph{large instances} composed of 100, 200 and 300 orders. Other four parameters, each taking three different values, are used to generated the instances. These parameters specify the layout of the warehouse, the storage policy, the trolley capacity, and the deadline distribution (see Table 3 of \cite{VanGils2019}). 
In the case of the small (resp.\ large) instances, for each set of the five parameters (the size of the orders plus the other four), 10 (resp.\ 30) instances are generated. This leads to a total of 2430 small (resp.\ 7290 large) instances. Thus, the benchmark of instances is made of  9720 instances. We report the results obtained with the procedure presented in Section~\ref{sec:algo} on the following subsets of instances:
\begin{itemize}
    \item small instances: we select only one replication for each set of parameters, thus 243 instances are considered;
    \item large instances: we select only instances with 100 orders and trolley capacity in $\{15,30\}$, for a total of $54$ instances.
\end{itemize}

We limit the small instances to one replication due to similar behaviour of the procedure over the other instances. For the large instances, the non-considered instances are clearly too large to be tackle\reven{d} with our approach. 

Some issues were found in the data sets of \cite{VanGils2019}. We contacted the authors and had kind and constructive exchanges with Kris Braekers who helped in correcting the problems.\footnote{The updated data-sets as well as the best solutions found are available on this page: \url{https://pagesperso.g-scop.grenoble-inp.fr/~cambazah/sequencing/}}
After the correction, one instance of the small as well as one of the large set turned out to be infeasible. We thus removed them for sake of comparison with the ILS proposed by \cite{VanGils2019}. Kris Braekers \reven{ran} the algorithm proposed in  \cite{VanGils2019} on the corrected benchmark. We thus compare here the performance of {\algo} against the new results obtained by the ILS. 
\reven{Because} we provide optimal results on all small instances, \reven{we can} assess the quality of the solutions obtained by the ILS.

The remainder of this section is organized as follows.
First, Section~\ref{sec:params} present\reven{s} the parameters used to run the algorithm.
In Sections~\ref{sec:resultsValidInequalities} and~\ref{sec:resultsTourCst}, we evaluate the interest of the proposed valid inequalities and strengthened tour constraints in the pricing, respectively.
Then, optimal results for the small instances are provided in Section~\ref{sec:small}.
The results of {\algo} on large instances are reported in Section~\ref{sec:large}.
Finally, in Section~\ref{sec:timelimit} \reven{we discuss} the impact of the time limit on the quality of the results of {\algo}.

\subsection{Parameters and configuration}
\label{sec:params}

The {\algo} uses a set of parameters list\reven{ed} here. 
The procedure that fills pool $\pool$ is run for a time $\tau_{pool}$. 
The resolution of $Pr(\mathcal{K}')$ has a time limit of $\tau_{Pr}$.
When {\algo} terminates, $(B(\mathcal K'))$ is solved with a time limit $\tau_l$ of the available remaining time to solve the problem.
\reven{If} $\tau_l$ is less than 600 seconds, $\tau_l$ is set to 600 seconds; this ensures a reasonable amount of time to solve $(B(\mathcal K'))$ for large instances.

In order to define the values of $\tau_{pool}$ and $\tau_{Pr}$, we \reven{first} run {\algo} on a subset of 9 large instances with different values for these two parameters.
These 9 instances (namely instances Id89, Id90, Id91, Id115, Id116, Id117, Id142, Id143 and Id144) \reven{have} 100 orders, trolley capacity $B=15$, and 24 aisles.
\reven{They are chosen because} they are challenging to solve (\reven{usually reaching} the time limit of 14~400 seconds), but the final optimality gaps are of acceptable quality.
After some preliminary experiments, we test values for $\tau_{pool}$ in $\{300; 600; 1~200\}$ seconds, and values for $\tau_{Pr}$ in $\{180; 300; 480\}$ seconds.
{\algo} is run with a time limit of 14~400 seconds for the 9 instances on each of the 9 configurations defined by the values of $\tau_{pool}$ and $\tau_{Pr}$.
The results are presented in Table~\ref{tableConfig}.
The first column report\reven{s} the time for $\tau_{pool}$.
Then three columns labeled \emph{\#TimeLimit} report the number of instances that reach the time limit, each column corresponding to a different value of $\tau_{Pr}$.
Finally, the last three columns labeled \emph{Avg. Gap(\%)} report the average optimality gap, each column corresponding to a different value of $\tau_{Pr}$.


\begin{table}[H]
\centering
{\footnotesize
\begin{tabular}{|c||r|r|r||r|r|r|}
\cline{2-7}
\multicolumn{1}{c||}{ }	&	\multicolumn{3}{c||}{\#TimeLimit}	&	\multicolumn{3}{c|}{Avg. Gap(\%)}	\\
\hline
$\tau_{pool}$ / $\tau_{Pr}$	&	180	&	300	&	480	&	180	&	300	&	480 \\
\hline
300	&	4	&	4	&	4	&	 0.57 	 & 	 0.58 	 & 	 0.54 	\\
600	&	4	&	4	&	5	&	 0.64 	 & 	 0.64 	 & 	 0.64 	\\
1~200	&	4	&	4	&	4	&	 0.52 	 & 	 0.53 	 & 	 0.52 	\\
\hline
\end{tabular}
}
\caption{Performance of CGH to solve 9 instances with 100 orders, trolley capacity $B=15$ and 24 aisles, with different parameters for $\tau_{pool}$ and $\tau_{Pr}$.}
\label{tableConfig}
\end{table}

Table~\ref{tableConfig} \reven{shows} that $\tau_{Pr}$ (the time limit to solve $Pr(\mathcal{K}')$) has almost no impact on the performance of {\algo}.
Regarding $\tau_{pool}$ (the time to run the procedure to fill the pool) it is clear that the best results are obtained for $\tau_{pool} = 1~200$ seconds.
Hence, it seems interesting to let this procedure run quite a long time to fill the pool with many \rev{good-quality} columns. So we set the values of $\tau_{pool}$ to 1~200 seconds, and $\tau_{Pr}$ to 180 seconds.


\subsection{Analysis on the contribution of the valid inequalities}
\label{sec:resultsValidInequalities}

In this section, we evaluate the potential of each family of inequalities presented in Sections~\ref{sec:ineqdff} and~\ref{sec:cutsjobprp}. To this end\rev{,} we solve $(M(\mathcal K))$, the linear relaxation of $(B(\mathcal{K}))$ on the small instances, with different configurations.
Since the instances are small, we can generate all columns in $\mathcal{K}$ and evaluate their time by calling the dynamic programming of~\cite{Cambazard18}.

Table \ref{tableSmallRes_cuts} provides details on the contribution of each family of cuts on the value of $(M(\mathcal K))$ the linear relaxation of $(B(\mathcal K))$. The first two columns labelled {\ttfamily Data sets} report the number of orders ($|\mathcal{O}|$) and the number of instances (\#Inst). The remaining columns, labelled {\ttfamily Average Root Gap (\%)} report the optimality gap at the root node of the \rev{branch-and-bound} tree for six different configurations: configuration labelled {\ttfamily None} solves $(\mathcal{M}(\mathcal K))$ without the use of any cut; configuration {\ttfamily All} uses all families of cuts; configuration {\ttfamily No cut}, with $cut \in \{MT, FS, R1, SC\} $ uses all the families of cuts but \reven{the} one indicated by parameter $cut$. More precisely, $cut = MT$ means Martello and Toth's cuts \eqref{cnst:l2bound} (Section \ref{sec:l2cuts}); $cut = FS$ means Fekete and Schepers's cuts \eqref{eq:vi01} (Section \ref{sec:FScuts}); $cut = R1$  means rank-1 cuts~\eqref{eq:rank1} (Section \ref{sec:rank1}); $cut = SC$ means the strengthened capacity cuts~\eqref{eq:sr} (Section \ref{sec:scccuts}). 
The {\ttfamily Average Root Gap (\%)} is computed as: $100 \frac{z^*_{B(\mathcal{K})} - z^*_{M
(\mathcal{K})}}{z^*_{M(\mathcal{K})}}$, where $z^*_{B(\mathcal{K})}$ is the optimal value of $B(\mathcal{K})$, and $z^*_{M(\mathcal{K})}$ is the optimal value of $M(\mathcal{K})$ after adding the corresponding violated cuts.

\begin{table}[H]
\centering
{\footnotesize
\begin{tabular}{|c|c||r|r|r|r|r|r|}
\hline
\multicolumn{2}{|c||}{Data sets} &\multicolumn{6}{c|}{Average Root Gap (\%)} \\
\hline
$|\mathcal O|$	&	\#Inst	&	None	&	All 
&	No MT	&	No FS	&	No R1	&	No SC	\\
\hline
6	&	81	&	0.35	&	0	&	0.14	&	0	&	0	&	0	\\
12	&	81	&	0.84	&	0.05	&	0.06	&	0.06	&	0.14	&	0.15	\\
18	&	80	&	1.00	&	0.08	&	0.14	&	0.08	&	0.27	&	0.33	\\
\hline
\end{tabular}
}
\caption{Average root optimality gap for different sets of valid inequalities.}
\label{tableSmallRes_cuts}
\end{table}

Overall the strengthened capacity cuts and the rank-1 cuts seem to be the most effective while the Fekete and Schepers's cuts have only a small impact on instances with 12 orders. It is also interesting to note that for instances with 6 orders, not adding the Martello and Toth's cuts yields a positive root gap while not adding other families of cuts has no impact. Moreover, we also note from our experiments that the strengthened capacity cut defined over the whole set of orders $\mathcal{O}$ is a very effective cut.
In conclusion, since all the families of inequalities help in solving the considered set of instances, we consider all of them in the final configuration of {\algo}.

\subsection{Interest of the strengthened tour constraints~\eqref{cst:cutpricing3}}
\label{sec:resultsTourCst}

In this section, we assess the interest of using the strengthened tour constraints in the MIP formulation of the pricing problem (see Section~\ref{sec:tourCst}).
On the 26 instances with 100 orders and a trolley capacity of 15, we have run {\algo} without considering the strengthened tour constraints, \ie by considering only the classical tour constraints~\eqref{cst:cutpricing} as already proposed in \cite{BriantCCCLO20}.
The time limit is set to 14~400 seconds.
Table~\ref{tableTourCst} reports overall results on this set of 26 instances.
Column {\ttfamily Strengh.TourCst.}\ indicates \reven{whether} {\algo} \reven{is} run with or without strengthened tour constraints.
The other columns report the number of instances that reach the time limit ({\ttfamily \#TimeLimit}), the average computation time in seconds to run all the instances ({\ttfamily Avg.\ Cpu(s)}), the number of instances solved to proven optimality ({\ttfamily \#Opt}), the number of instances with a final gap that is lower than 1\% ({\ttfamily \#Gap$\leq$1\%}) and the average optimality gap ({\ttfamily Avg.\ Gap(\%)}).

\begin{table}[H]
\centering
{\footnotesize
\begin{tabular}{|c||r|r|r|r|r|}
\hline
Strengh.TourCst.	&	\#TimeLimit	&	Avg. Cpu(s)	&	\#Opt & \#Gap$\leq$1\%	&	Avg. Gap(\%)	\\
\hline
yes	&	9	&	8 399	&	4	&	22	&	0.74	\\
no	&	17	&	12 231	&	3	&	12	&	3.73	\\
\hline
\end{tabular}
}
\caption{Overall comparison of CGH on data sets of 100 orders and batch capacity 15, with and without strengthened tour constraints in the pricing.}
\label{tableTourCst}
\end{table}

Table~\ref{tableTourCst} \reven{clearly shows} that adding the strengthened tour constraints in the MIP formulation of the pricing problem \reven{leads to shorter} computation time\reven{s} and optimality gaps of higher quality. \reven{This} reflect\reven{s} the theoretical dominance of the strengthened tour constraints~\eqref{cst:cutpricing3} over the tour constraints~\eqref{cst:cutpricing} from a computational point of view. We th\reven{us} consider them in the final configuration of {\algo}.

\subsection{Exact algorithm using $(B({\cal K}))$ for small size instances}
\label{sec:small}

On the small instances a complete enumeration of the feasible {\tour}s is possible in a reasonable amount of time. 
The duration of each {\tour} is computed with the dynamic programming of~\cite{Cambazard18}.
Therefore formulation $(B({\mathcal K}))$, strengthened with all the valid inequalities presented in Sections~\ref{sec:ineqdff} and~\ref{sec:cutsjobprp}, can be solved as a compact integer programming formulation.\footnote{Note that when ignoring the sequencing part of the problem, it boils down to the batching and picker routing problem. All small instances can also be optimally solved in this case. 
} 

The results on small instances are shown in Table \ref{tableSmallRes}. The first two columns labelled {\ttfamily Data sets} report the number of orders ($|\mathcal{O}|$) and the number of instances ({\ttfamily \#Inst}). The next five columns labelled {\ttfamily Cpu\ (s)} report, in seconds, the following computational times: the average ({\ttfamily avg}), the minimum ({\ttfamily min}) and the maximum ({\ttfamily max}) computational time to solve one instance with the corresponding characteristics, the average time to generate all the {\tour}s ({\ttfamily InitGen}) and the average time to solve the MIP formulation ({\ttfamily MIP}). 

All instances can be solved to optimality within the given time limit. The maximum time needed is 410.2 seconds. The generation of all columns is the most time-consuming part and takes \rev{on} average 40.0 seconds for instances with 18 orders. The resulting MIP formulation is often solved relatively quickly: it takes 4.9 seconds \rev{on} average for instances with 18 orders. The remaining time (not detailed in the table) is used to separate the family of cuts that we consider. \reven{These results} considerably outperform the exact approach proposed in \cite{VanGils2019} which fails to solve to proven optimality 40.6\% of the small instances within a time limit of 4 hours.
\reven{Incidentally, } having optimal solutions for all small instances now enables to assess the quality of the ILS proposed in \cite{VanGils2019}.


\begin{table}[H]
\centering
{\footnotesize
\begin{tabular}{|c|c||r|r|r||r|r|}
\hline
\multicolumn{2}{|c||}{Data sets} &\multicolumn{5}{c|}{Cpu (s)}   \\
\hline
|$\mathcal O$|& \#Inst  & avg & min & max & InitGen & MIP \\
\hline
6	&	81	&	0	&	0.0	&	0.1	&	0	&	0	\\
12	&	81	&	1.5	&	0.0	&	6.5	&	1.4	&	0	\\
18	&	80	&	49.4	&	0.0	&	410.2	&	40.0	&	4.9	\\
\hline
\end{tabular}
}
\caption{Overall results on the small instances. }
\label{tableSmallRes}
\end{table}



To conclude, formulation $(B({\cal K}))$ can solve to optimality all instances with up to 18 orders. This is thanks to the fact that all the {\tour}s can be generated beforehand. 

\subsection{Evaluation of {\algo} on large size instances}
\label{sec:large}

We turn our attention to the ability of {\algo} to provide \rev{high-quality} lower and upper bounds on the large instances.
A time limit of 14~400 seconds (4 hours) is imposed \reven{-- but} the final resolution of $B(\mathcal{K'})$ is not included in the time limit, hence the final cpu time may exceed 4 hours.

Tables~\ref{tableLargeRes15} and~\ref{tableLargeRes30} report the results obtained on instances of $100$ orders and a trolley capacity $\capa$ equal to 15 and 30. 

The first five columns labelled {\ttfamily Data set} report the identification of the instance ({\ttfamily Id}), the number of pickers ({$|\mathcal P|$}), the number of aisles in the warehouse ({\ttfamily \#Aisles}), the storage location policy  \reven{used when generating} the instance ({\ttfamily Loc}) and the deadline distribution ({\ttfamily DD}). Values reported in  columns {\ttfamily Loc} and {\ttfamily DD} are in $\{1, 2, 3\}$ with the following meaning: 1=\lq Random\rq, 2=\lq Within aisle\rq, 3=\lq Across aisle\rq\ and 1=\lq Uniform\rq, 2=\lq Triangular progressive\rq, 3=\lq Triangular degressive\rq, respectively. 


The sixth column ({\ttfamily UB}) reports the value of the best solution found by the ILS algorithm of \cite{VanGils2019}. 
The last six columns labelled {\ttfamily \algo} \reven{detail} the performance of {\algo}. 
In particular column {\ttfamily UB init.}\ reports the upper bound found by the initial heuristic (see Section~\ref{sec:initmaster}), column {\ttfamily UB} reports the best upper bound found by {\algo}, column {\ttfamily Cpu~(s)} indicates the computation time in seconds, column {\ttfamily Gap~(\%)} reports the optimality gap computed as: $100 \frac{z^{UB} - z^{LB}}{z^{LB}}$, where $z^{UB}$ and $z^{LB}$ respectively denote the best upper and lower bounds obtained at the end of CGH.
The optimality gap is reported in italic when it is zero, to emphasize that the optimality has been proven on the respective instance. Finally, columns {\ttfamily \#Cols} and {\ttfamily \#It} report the number of columns in $\mathcal K'$ at the end of the procedure and the number of calls to the pricing problem (see line~\ref{lpr1} of Algorithm~\ref{alg:genColAlgo}), respectively.
Note that when some negative reduced cost columns are found in the pool, this is not counted as an extra iteration.
%


\begin{table}[ht]
\centering
{\footnotesize
\begin{tabular}{|c|c|c|c|c||r||r|r|r|r|r|r|r|}
\hline
\multicolumn{5}{|c||}{Data set} &\multicolumn{1}{|c||}{ILS} &\multicolumn{6}{c|}{\algo} \\
\hline
Id	&	|P|	&	\#Aisles	&	Loc	&	DD	&	UB	& UB init. &	UB	&	Cpu (s)	& Gap (\%)	&	\#Cols	&	\#It	\\
\hline
7	&	4	&	12	&	1	&	1	&	 360 314 	&	 359 932 	&	 \textbf{359 348} 	&	 3 508 	&	0.02	&	 5 303 	&	35	\\
8	&	3	&	12	&	1	&	2	&	 308 226 	&	 307 526 	&	 \textbf{307 070} 	&	 7 685 	&	0.04	&	 11 251 	&	38	\\
9	&	4	&	12	&	1	&	3	&	 337 286 	&	 334 176 	&	 \textbf{332 952} 	&	 4 116 	&	0.01	&	 5 792 	&	28	\\
34	&	3	&	12	&	2	&	1	&	 285 768 	&	 282 386 	&	 \textbf{281 710} 	&	 3 257 	&	0.04	&	 4 775 	&	28	\\
35	&	3	&	12	&	2	&	2	&	 279 886 	&	 279 218 	&	 \textbf{278 818} 	&	 1 586 	&	0	&	 4 605 	&	17	\\
36	&	4	&	12	&	2	&	3	&	 265 444 	&	 265 172 	&	 \textbf{264 646} 	&	 3 421 	&	0.05	&	 7 088 	&	45	\\
61	&	3	&	12	&	3	&	1	&	 292 240 	&	 291 174 	&	 \textbf{290 644} 	&	 9 267 	&	0.01	&	 10 412 	&	36	\\
62	&	4	&	12	&	3	&	2	&	 337 164 	&	 337 060 	&	 \textbf{336 990} 	&	 2 039 	&	0	&	 3 763 	&	13	\\
63	&	4	&	12	&	3	&	3	&	 294 490 	&	 291 382 	&	 \textbf{290 808} 	&	 6 814 	&	0.01	&	 7 193 	&	31	\\
88	&	5	&	24	&	1	&	1	&	 510 786 	&	 505 468 	&	 \textbf{504 440} 	&	 14 407 	&	0.08	&	 9 011 	&	24	\\
89	&	5	&	24	&	1	&	2	&	 453 518 	&	 450 730 	&	 \textbf{448 880} 	&	 14 403 	&	0.42	&	 7 914 	&	18	\\
90	&	6	&	24	&	1	&	3	&	 533 726 	&	 533 164 	&	 \textbf{531 770} 	&	 14 414 	&	2.32	&	 8 537 	&	17	\\
115	&	5	&	24	&	2	&	1	&	 440 786 	&	 440 710 	&	 \textbf{437 212} 	&	 2 800 	&	0.01	&	 4 839 	&	10	\\
116	&	5	&	24	&	2	&	2	&	 487 356 	&	 485 936 	&	 \textbf{482 316} 	&	 2 727 	&	0.02	&	 4 092 	&	9	\\
117	&	6	&	24	&	2	&	3	&	 426 670 	&	 427 634 	&	 \textbf{426 282} 	&	 4 509 	&	0.07	&	 4 967 	&	20	\\
142	&	5	&	24	&	3	&	1	&	 452 266 	&	 452 142 	&	 \textbf{446 276} 	&	 13 463 	&	0.02	&	 7 867 	&	28	\\
143	&	5	&	24	&	3	&	2	&	 449 182 	&	 446 764 	&	 \textbf{443 084} 	&	 5 111 	&	0	&	 5 408 	&	12	\\
144	&	6	&	24	&	3	&	3	&	 471 632 	&	 470 744 	&	 \textbf{470 066} 	&	 14 424 	&	1.79	&	 8 228 	&	19	\\
169	&	6	&	36	&	1	&	1	&	 708 648 	&	 701 362 	&	 \textbf{694 036} 	&	 15 001 	&	0.55	&	 7 700 	&	22	\\
170	&	6	&	36	&	1	&	2	&	 583 298 	&	 581 488 	&	 \textbf{579 718} 	&	 14 422 	&	11.79	&	 10 060 	&	5	\\
196	&	6	&	36	&	2	&	1	&	 559 908 	&	 557 294 	&	 \textbf{553 256} 	&	 14 415 	&	0.05	&	 7 313 	&	12	\\
197	&	6	&	36	&	2	&	2	&	 652 180 	&	 650 868 	&	 \textbf{650 470} 	&	 4 315 	&	0.01	&	 4 829 	&	8	\\
198	&	7	&	36	&	2	&	3	&	 547 630 	&	 548 534 	&	 \textbf{541 890} 	&	 10 908 	&	0.02	&	 8 914 	&	13	\\
223	&	6	&	36	&	3	&	1	&	 525 258 	&	 518 040 	&	 \textbf{517 204} 	&	 14 404 	&	1.89	&	 8 307 	&	12	\\
224	&	6	&	36	&	3	&	2	&	 614 394 	&	 614 510 	&	 \textbf{610 220} 	&	 14 345 	&	0.03	&	 7 414 	&	20	\\
225	&	7	&	36	&	3	&	3	&	 654 416 	&	 653 408 	&	 \textbf{649 296} 	&	 2 628 	&	0	&	 4 043 	&	2	\\
\hline
\end{tabular}
}
\caption{Detailed results on data sets of 100 orders ($|\mathcal{O}| = 100$), and trolley capacity $B =15$.}
\label{tableLargeRes15}
\end{table}

Table~\ref{tableLargeRes15} \reven{supports} the following observations.
\begin{itemize}
    \item 17 out of the 26 instances \rev{provided an optimal solution of} the linear relaxation of $B(\mathcal{K})$ before the time limit of 4 hours.
    \item 22 out of the 26 instances are solved to near optimality, \reven{with an} optimality gap at the end of the computation lower than 1\%. Four instances, \reven{namely} Id35, Id62, Id143 and Id225, are solved to optimality. We can claim that the lower bounds provided by $(LM(\mathcal{K}))$ are of very high quality even on large size data sets. 
    \item The lower bound computed by {\algo} allows to assert the quality of ILS on this benchmark.
    \reven{ILS is overall very good}: for the 20 \rev{instances} with less than 0.1\% of optimality gap, ILS provides solutions with an average deviation of 0.82\% from the lower bounds, and a maximal deviation of 1.48\%.
    \item The initial heuristic provides good solutions, often better than the ones provided by ILS, but the time for the heuristic is larger than the one \reven{of} ILS since the idea is to generate an interesting pool of columns. For three instance\reven{s}, the solution of the initial heuristic is wors\reven{e} than the one of ILS. 
    In all instances the final upper bound of {\algo} is better than the initial upper bound.
    \item In comparison with the results of ILS, all upper bounds are improved by {\algo} but with a significantly longer time since {\algo} also provides lower bounds (ILS is run with a limited number of iterations and its runtime is 22 seconds \rev{on} average with a maximum of 70 seconds).
    \reven{On} average, the upper bounds reported by {\algo} decrease by 0.83\% the ones of ILS, with a maximum of 2.06\% for instance Id169.
    \item The number of calls to the pricing problem is rather small for the given time limit. 
    Note that negative reduced cost columns are found in the pool 75 times \reven{on} average (see Line~\ref{line:checkPool} of Algorithm~\ref{alg:genColAlgo}).
    These are not counted as iterations since there is no call to the pricing algorithm.
    Moreover, more than 95\% of the time is spent to solve the pricing problem, hence the solving of the pricing problem is the bottleneck of the proposed approach.
\end{itemize}

Table~\ref{tableLargeRes30} reports results obtained on instances with trolley capacity of 30.

\begin{table}[ht]
\centering
{\footnotesize
\begin{tabular}{|c|c|c|c|c||r||r|r|r|r|r|r|r|}
\hline
\multicolumn{5}{|c||}{Data set} &\multicolumn{1}{|c||}{ILS} &\multicolumn{6}{c|}{\algo} \\
\hline
Id	&	|P|	&	\#Aisles	&	Loc	&	DD	& UB  &	UB init.	&	UB	&	Cpu (s)	& Gap (\%)	&	\#Cols	&	\#It	\\
\hline
16	&	4	&	12	&	1	&	1	&	 305 840 	&	 305 766 	&	 \textbf{304 360} 	&	 15 039 	&	2.67	&	 13 724 	&	37	\\
17	&	3	&	12	&	1	&	2	&	 268 896 	&	 267 480 	&	 \textbf{266 792} 	&	 14 554 	&	5.54	&	 16 242 	&	19	\\
18	&	4	&	12	&	1	&	3	&	 256 420 	&	 256 378 	&	 \textbf{255 388} 	&	 14 995 	&	12.15	&	 17 460 	&	19	\\
43	&	3	&	12	&	2	&	1	&	 219 278 	&	 219 116 	&	 \textbf{218 306} 	&	 15 001 	&	11	&	 17 111 	&	9	\\
44	&	3	&	12	&	2	&	2	&	 203 742 	&	 203 482 	&	 \textbf{203 022} 	&	 15 117 	&	6.06	&	 22 115 	&	28	\\
45	&	4	&	12	&	2	&	3	&	 201 350 	&	 202 194 	&	 \textbf{200 858} 	&	 14 628 	&	1.59	&	 23 992 	&	84	\\
70	&	3	&	12	&	3	&	1	&	 232 510 	&	 231 860 	&	 \textbf{230 976} 	&	 14 495 	&	9.68	&	 17 633 	&	39	\\
71	&	3	&	12	&	3	&	2	&	 \textbf{228 188} 	&	 228 572 	&	 228 342 	&	 14 502 	&	7.64	&	 14 716 	&	19	\\
72	&	4	&	12	&	3	&	3	&	 239 470 	&	 239 220 	&	 \textbf{238 016} 	&	 14 572 	&	12.89	&	 15 114 	&	5	\\
97	&	4	&	24	&	1	&	1	&	 343 062 	&	 341 960 	&	 \textbf{339 116} 	&	 15 010 	&	96.35	&	 18 975 	&	3	\\
98	&	5	&	24	&	1	&	2	&	 458 204 	&	 456 140 	&	 \textbf{454 030} 	&	 14 804 	&	34.76	&	 15 362 	&	9	\\
99	&	5	&	24	&	1	&	3	&	 417 336 	&	 414 356 	&	 \textbf{413 566} 	&	 15 008 	&	84.12	&	 19 666 	&	23	\\
124	&	4	&	24	&	2	&	1	&	 \textbf{322 434} 	&	 322 608 	&	 322 608 	&	 15 016 	&	20.67	&	 16 276 	&	8	\\
125	&	4	&	24	&	2	&	2	&	 321 718 	&	 319 336 	&	 \textbf{318 336} 	&	 15 004 	&	16.46	&	 15 744 	&	6	\\
126	&	5	&	24	&	2	&	3	&	 349 442 	&	 345 970 	&	 \textbf{345 340} 	&	 14 855 	&	29.08	&	 19 766 	&	13	\\
151	&	4	&	24	&	3	&	1	&	 358 618 	&	 357 596 	&	 \textbf{357 596} 	&	 15 025 	&	53.95	&	 16 389 	&	5	\\
152	&	4	&	24	&	3	&	2	&	 361 010 	&	 359 750 	&	 \textbf{359 054} 	&	 15 028 	&	28.61	&	 16 945 	&	6	\\
153	&	5	&	24	&	3	&	3	&	 383 112 	&	 378 886 	&	 \textbf{377 676} 	&	 14 437 	&	14.34	&	 13 205 	&	14	\\
178	&	5	&	36	&	1	&	1	&	 480 888 	&	 476 112 	&	 \textbf{476 112} 	&	 15 010 	&	129.74	&	 17 898 	&	5	\\
179	&	5	&	36	&	1	&	2	&	 531 724 	&	 523 032 	&	 \textbf{521 038} 	&	 15 041 	&	101.25	&	 16 001 	&	4	\\
180	&	6	&	36	&	1	&	3	&	 600 422 	&	 576 516 	&	 \textbf{576 516} 	&	 15 022 	&	76.5	&	 14 579 	&	3	\\
205	&	5	&	36	&	2	&	1	&	 451 856 	&	 449 758 	&	 \textbf{445 436} 	&	 14 931 	&	65.69	&	 15 814 	&	3	\\
206	&	6	&	36	&	2	&	2	&	 502 722 	&	 504 404 	&	 \textbf{498 140} 	&	 15 042 	&	40.59	&	 14 589 	&	4	\\
207	&	7	&	36	&	2	&	3	&	 503 854 	&	 500 510 	&	 \textbf{497 284} 	&	 15 035 	&	45.61	&	 14 290 	&	11	\\
232	&	6	&	36	&	3	&	1	&	 523 364 	&	 522 192 	&	 \textbf{521 036} 	&	 15 005 	&	79.16	&	 14 745 	&	4	\\
233	&	6	&	36	&	3	&	2	&	 497 086 	&	 497 472 	&	 \textbf{495 788} 	&	 15 028 	&	38.16	&	 13 921 	&	10	\\
234	&	6	&	36	&	3	&	3	&	 450 494 	&	 446 924 	&	 \textbf{445 040} 	&	 15 013 	&	71.38	&	 16 399 	&	3	\\
\hline
\end{tabular}
}
\caption{Detailed results on data sets of 100 orders ($|\mathcal{O}| = 100$), and trolley capacity $B = 30$.}
\label{tableLargeRes30}
\end{table}

The results in Table~\ref{tableLargeRes30} \reven{point to} the following observations.
\begin{itemize}
    \item All instances ran up to the time limit of 4 hours.
    \item Two instances (Id16 and Id45) \reven{yield} a gap of less than 3\%, and six instances \reven{yield} a gap of less than 10\%. For the remaining instances, the optimality gaps are rather large.
    Hence, as expected, instances with $B = 30$ are much more difficult to be solved with {\algo}.
    \item However, {\algo} is able to improve the value of the upper bound with respect to the ILS for all instances except two (Id71 and Id124). 
    Note that for these two instance, {\algo} report\reven{s} an upper bound very close from the ones \reven{of} ILS.
    For the other instances, the upper bounds of {\algo} are \rev{on} average 0.95\% less than the ones of ILS, and with some significant improvement for instances Id179 and Id180.
\end{itemize}

As a conclusion, CGH provide\reven{s} overall very good results on large size instance with a small trolley capacity ($B=15$).
However, the performance of CGH deteriorate when enlarging the capacity of the trolley to $B=30$.
In such a case the pricing problem is more difficult to solve, \reven{hence} obtain\reven{ing} a \rev{good-quality} lower bound \reven{becomes difficult}.

\subsection{Impact of the time limit to the performance of CGH}
\label{sec:timelimit}

In this section, we evaluate the impact of the time limit on the quality of the solutions reported by {\algo}.
On the 26 instances with 100 orders and a trolley capacity $B=15$, we have run {\algo} with different time limits: 1~800 seconds, 3~600 seconds, 7~200 seconds and 14~400 seconds.
Note that the time to fill the pool $\tau_{pool} = 1~200$ seconds is included in the solving time.
Table~\ref{tableTimeLimit} reports overall results on this set of 26 instances.  
Column \texttt{Time Limit(s)} indicates the time limit to run {\algo} in seconds.
The other columns report the number of instances that reach the time limit (\texttt{\#TimeLimit}), the number of instances solved to proven optimality (\texttt{\#Opt}), the number of instances with a final gap that is lower than 1\% (\texttt{\#Gap$\leq$1\%}) and the average optimality gap (\texttt{Avg.\ Gap(\%)}).

\begin{table}[H]
\centering
{\footnotesize
\begin{tabular}{|c||r|r|r|r|}
\hline
Time Limit(s)	&	\#TimeLimit	&	\#Opt	&	\#Gap$\leq$1\%	&	Avg. Gap(\%)	\\
\hline
1 800	&	25	&	1	&	3	&	14.64	\\
3 600	&	18	&	3	&	12	&	3.98	\\
7 200	&	12	&	4	&	16	&	1.79	\\
14 400	&	9	&	4	&	22	&	0.74	\\

\hline
\end{tabular}
}
\caption{Impact of the time limit \reven{on} the performance of CGH for 26 instances with 100 orders and trolley capacity $B=15$.}
\label{tableTimeLimit}
\end{table}

Table~\ref{tableTimeLimit} shows that the results are already of acceptable quality with\reven{in} a one hour time limit.

\subsection{Evaluation of the benefit of integrating the decisions}

\label{sec:benefitIntegrated}

In this section, we aim to evaluate the advantage of integrating sequencing decisions with order batching and picker routing decisions when solving the \pbname{}.
We compare the integrated approach with a wave picking system, \ie a very classical approach where sequencing decisions are not fully integrated with order batching and picker routing decisions. 
Recall that Figure~\ref{fig:sequencing} shows an example of a solution obtained in a wave picking system compared to a solution of the \pbname{}.
In order to simulate a wave picking system, we consider a restriction of the \pbname{} where a batch $k$ is considered valid if all collected orders $o\in\mathcal{O}_k$ have the same deadline, \ie $\bar{d}_k = \min_{o \in\mathcal{O}_k}{\bar{d}_o} = \max_{o \in\mathcal{O}_k}{\bar{d}_o}$.

Table~\ref{tableWavePicking} reports a comparison of the solutions obtained when solving the \pbname{} with or without the wave picking restriction on the set of 242 small instances.
The first two columns, entitled {\ttfamily Data sets}, indicate the number of orders ($|\mathcal{O}|$) and the number of instances ({\ttfamily \#Inst}).
The next two columns report the number of feasible solutions ({\ttfamily \#Feasible}) and the average optimal values ({\ttfamily Avg.\ Opt}) when the \pbname{} is solved without wave picking restriction.
The subsequent two columns report the same information when the wave picking restriction is considered.
In the wave picking case, the average optimal values are only computed for instances reporting a feasible solution.
Finally, column \texttt{Deviation~(\%)} reports the percentage deviation of the average total duration when taking the decisions in an integrated way when solving the \pbname.

\begin{table}[H]
\centering
{\footnotesize
\begin{tabular}{|c|c||r|r|r|r|r|}
\hline
\multicolumn{2}{|c||}{Data sets}	&	\multicolumn{2}{c|}{\pbname{}} & \multicolumn{2}{c|}{Wave picking} & \multirow{2}{*}{Deviation (\%)}	\\
\cline{1-6}
$|\mathcal{O}|$	&	\#Inst	&	\#Feasible	&	Avg. Opt	&	\#Feasible	&	Avg. Opt	&		\\
\hline
6	&	81	&	81	&	 75~981 	&	80	&	 92~167 	&	-17.6 	\\
12	&	81	&	81	&	 90~017 	&	79	&	 120~988 	&	-25.6 	\\
18	&	80	&	80	&	 101~668 	&	79	&	 138~662 	&	-26.7 	\\
\hline
\end{tabular}
}
\caption{Comparison of optimal values on the small size instances for the \pbname{} and a wave picking approach.}\label{tableWavePicking}
\end{table}

Table~\ref{tableWavePicking} shows that adopting a wave picking system instead of an integrated solving approach leads to non feasible solutions in a few cases only (4 out of 242).
However, the main interest of taking the decisions in an integrated way is the very significant decrease of the total time.

\section{Conclusion}\label{sec:conclusion}

\reven{This paper shows} that the Joint Order Batching, Picker Routing and Sequencing Problem with Deadlines (\pbname) is better captured by formulating it as a bin packing problem rather than a scheduling problem. We believe that the heuristic approaches proposed in \cite{VanGils2019} can be improved based on this analysis. This observation lead\reven{s} to an algorithm with performance guarantee (\ie able to provide valid lower and upper bounds), for a very complex integrated logistic problem. \rev{Experimental} results show that the proposed approach correctly scales \reven{up} on instances of reasonable sizes.

A second contribution is the design of valid inequalities for the master problem as well as for the pricing problem. The first family of cutting planes takes advantage of the bin packing analysis and relies on Dual-Feasible Functions. Regarding the pricing problem, the inequalities proposed in \cite{BriantCCCLO20} (namely the tour constraints) are also strengthened. This allows to efficiently solve the pricing problem \reven{perceived} as the  bottleneck of the approach due to the hardness of the picker routing problem. As a result, this paper shows that the column generation based heuristic (\algo) proposed for the Joint Order Batching and Picker Routing Problem (JOBPRP) in \cite{BriantCCCLO20} can be extended to consider deadlines associated with \reven{the} orders to be prepared.

The proposed algorithm is able to solve optimally all the small instances generated by \cite{VanGils2019} and involving up to 18 orders. It provides very tight intervals of the optimal values (and improved upper bounds) for some of the large instances of 100 orders. It therefore contributes in asserting the quality of the heuristic techniques often proposed in this area of research due to the hardness of \reven{such} integrated problems.

In terms of future research, a first direction is to improve the solving of the pricing problem when dealing with large size instances, \ie many orders and high trolley capacity. The solving of the pricing problem is clearly the bottleneck to address such large instances, hence future improvements likely lie in this component of \algo .

The next avenue of research on the topic is to integrate additional constraints related to sequencing (release dates, breaks for pickers, dynamic setting such as \cite{DHaen2022}) or additional decision levels such as storage location assignment.
But we believe that a simpler and better algorithmic framework for the pricing \reven{is needed} before extending  the problem scope further.

Future research may benefit from the ideas presented in this paper to develop efficient heuristics that take advantage of the bin-packing structure of the problem.
A fast heuristic version of \algo{} could also be considered for implementation in an industrial context.

\bibliographystyle{apalike-ejor}
\bibliography{biblio}
\newpage

\section{Appendix}

\subsection{Proof of Proposition \ref{prop:L2cuts}}
\label{proof:L2cuts}

\paragraph{Proposition \ref{prop:L2cuts}}

{\em The following constraints}
\begin{equation*}
\sum_{k \in \mathcal{K}_1(\bar{d},q)} \bar{d} \, \rho_{k} + \sum_{k \in \mathcal{K}_2(\bar{d},q)} t_k \, \rho_{k} \leq \bar{d} \times |\mathcal{P}| \;\;\;\; \forall
\bar{d} \in \mathcal{D}, q \in \mathcal{Q}(\bar{d})
 \qquad \eqref{cnst:l2bound} 
\end{equation*}
{\em are valid inequalities for $(M(\mathcal{K}))$.}

\begin{proof}
Consider the function $f_0^{\lambda}$ from the $L_2$ lower bound proposed by \cite{Martello1990}. Let $\lambda \in \left[0;\frac{1}{2}\right]$:
\begin{align*}
f_0^{\lambda} : \left[ 0,1\right] &\rightarrow \left[0,1\right] \\
x & \mapsto \left\{
\begin{array}{l}
  1, \;\; \text{if } x > 1-\lambda \\
  x, \;\; \text{if } \lambda \leq x \leq 1-\lambda \\
  0, \;\; \text{if } x < \lambda
\end{array}
\right.
\end{align*}

By applying Proposition \ref{prop:DFF} to constraints~(\ref{cnst2BP}) with DFF $f_0^{\lambda}$, the following valid inequalities are obtained:
\begin{equation*}
\sum_{k \in \mathcal{K}(\bar{d})}{f_0^{\lambda}\left(\frac{t_k}{\bar{d}}\right) \rho_{kp}} \leq 1, \;\;\; \forall \bar{d} \in \mathcal{D}, p \in \mathcal{P}
\end{equation*}

This can be written as:
\begin{equation*}
\sum_{k \in \mathcal{K}(\bar{d}) ~|~ \frac{t_k}{\bar{d}} < \lambda}{0 \rho_{kp}} + \sum_{k \in \mathcal{K}(\bar{d}) ~|~ \lambda \leq \frac{t_k}{\bar{d}} \leq 1- \lambda}{\frac{t_k}{\bar{d}} \rho_{kp}} + \sum_{k \in \mathcal{K}(\bar{d}) ~|~ \frac{t_k}{\bar{d}} > 1- \lambda}{1 \rho_{kp}} \leq 1, \;\;\; \forall \bar{d} \in \mathcal{D}, p \in \mathcal{P}
\end{equation*}

By multiplying by $\bar{d}$, and by setting $\lambda = \frac{q}{\bar{d}}$, we obtain the following:
\begin{equation*}
\sum_{k \in \mathcal{K}(\bar{d}) ~|~ q \leq t_k \leq \bar{d} -q}{t_k \rho_{kp}} + \sum_{k \in \mathcal{K}(\bar{d}) ~|~ t_k > \bar{d}- q}{\bar{d} \rho_{kp}} \leq \bar{d}, \;\;\; \forall \bar{d} \in \mathcal{D}, p \in \mathcal{P}
\end{equation*}

By summing over all $p \in \mathcal{P}$ and taking into consideration \eqref{cnst:rhokp}, we get exactly constraints~(\ref{cnst:l2bound}).
\end{proof}

\subsection{Proof of Proposition \ref{prop:FeketeCuts}}
\label{proof:FScuts}

\paragraph{Proposition~\ref{prop:FeketeCuts}} {\em The following inequalities}
\begin{equation*}
\sum_{i=1}^{q-1} \sum_{k\in \mathcal{K}(\bar{d},i, q)} i \, \rho_k \leq (q-1) \times |\mathcal{P}|\quad\forall \bar{d}\in\mathcal{D}, q\in\{2,\dots,\bar{d}\}\quad \eqref{eq:vi01} 
\end{equation*}
{\em are valid for $(M(\mathcal{K}))$.}

\begin{proof}
\cite{Fekete2001} proposed a DFF denoted by $f_{FS,1}^{\lambda}$ (see \cite{Clautiaux2010}). Given a $\lambda \in \N \setminus \{0\}$, the function is the following:
\begin{align*}
f_{FS,1}^{\lambda} : \left[ 0,1\right] &\rightarrow \left[0,1\right] \\
x & \mapsto \left\{
\begin{array}{ll}
  x, & \text{if }  x (\lambda +1) \in \Z \\
  \left\lfloor (\lambda+1) x \right\rfloor \frac{1}{\lambda}, & \text{otherwise}
\end{array}
\right.
\end{align*}

We consider a weaker version of $f_{FS,1}^{\lambda}$ referred to as $g^{\lambda}$ which handles differently the case $x (\lambda +1) \in \Z$:

\begin{align*}
g^{\lambda}(x) = & \left\{
\begin{array}{ll}
    0  & \text{if }  x = 0\\
   \frac{x(\lambda + 1) - 1}{\lambda} & \text{if }  x (\lambda +1) \in \Z\setminus\{0\} \\
  \left\lfloor (\lambda+1) x \right\rfloor \frac{1}{\lambda}, & \text{otherwise}
\end{array}
\right.
\end{align*}
Note that $g^{\lambda}(x) \leq f_{FS,1}^{\lambda}(x)$ for any $x \in [0,1]$ so that $g^{\lambda}$ is a non-maximal DFF.
By applying Proposition \ref{prop:DFF} to constraints~(\ref{cnst2BP}) with $g^{\lambda}(x)$, the following valid inequalities are obtained:

\begin{equation*}
\sum_{k \in \mathcal{K}(\bar{d}) ~|~ \frac{t_k}{\bar{d}}(\lambda+1) \in \Z}{\frac{\frac{t_k}{\bar{d}}(\lambda+1) -1}{\lambda} \rho_{kp}} + \sum_{k \in \mathcal{K}(\bar{d}) ~|~ \frac{t_k}{\bar{d}}(\lambda+1) \notin \Z}{\left\lfloor (\lambda+1) \frac{t_k}{\bar{d}} \right\rfloor \frac{1}{\lambda} \rho_{kp}} \leq 1, \;\;\; \forall \bar{d} \in \mathcal{D}, p \in \mathcal{P}
\end{equation*}

Let us now set $\lambda = q-1$, and $\frac{t_k}{\bar{d}}(\lambda+1) = j$ and multiply by $q-1$. We obtain:

\begin{equation*}
\sum_{k \in \mathcal{K}(\bar{d}) ~|~ t_k = j \frac{\bar{d}}{q}, j \in \Z}{(j-1) \rho_{kp}} + \sum_{k \in \mathcal{K}(\bar{d}) ~|~ t_k = j \frac{\bar{d}}{q}, j \notin \Z}{\left\lfloor j \right\rfloor \rho_{kp}} \leq q-1, \;\;\; \forall \bar{d} \in \mathcal{D}, p \in \mathcal{P}
\end{equation*}

Consider $i = \left\lfloor j \right\rfloor$, we can replace $t_k = j \frac{\bar{d}}{q}, j \notin \Z$ by $i\frac{\bar{d}}{q} < t_k < (i+1) \frac{\bar{d}}{q}, i \in \Z$.
This gives the following constraints:
\begin{equation*}
\sum_{k \in \mathcal{K}(\bar{d}) ~|~ t_k = i \frac{\bar{d}}{q}, i \in \Z}{(i-1) \rho_{kp}} + \sum_{k \in \mathcal{K}(\bar{d}) ~|~ i\frac{\bar{d}}{q} < t_k < (i+1) \frac{\bar{d}}{q}, i \in \Z}{ i \rho_{kp}} \leq q-1, \;\;\; \forall \bar{d} \in \mathcal{D}, p \in \mathcal{P}
\end{equation*}

Thus the following inequality is valid:

\begin{equation*}
 \sum_{k \in \mathcal{K}(\bar{d}) ~|~ i\frac{\bar{d}}{q} < t_k \leq (i+1) \frac{\bar{d}}{q}, i \in \Z}{ i \rho_{kp}} \leq q-1, \;\;\; \forall \bar{d} \in \mathcal{D}, p \in \mathcal{P}
\qquad \constraintlabel{eq:fs1}
\end{equation*}

By summing (\ref{eq:fs1}) over all $p \in \mathcal{P}$ and taking into consideration \eqref{cnst:rhokp}, we get exactly constraints~(\ref{eq:vi01}).

\end{proof}

\subsection{Proof of Proposition \ref{prop:validPricingCuts}}
\label{proof:validPricingCuts}

\paragraph{Proposition \ref{prop:validPricingCuts}} {\em If constraint~\eqref{eq:eq2b} is valid and hypothesis~\eqref{eq:eq2bhyp} holds, then the following constraint}
\begin{equation*}
    b  + \sum_{o \in \mathcal{O}} a_o \,e_o \;+\; \Delta_k \left( \sum_{o \in\mathcal{O}_k} e_o - |k| + 1 \right) \leq t
    \qquad \eqref{cst:cutpricing3}
\end{equation*}
{\em is valid for $Pr(\mathcal{K}')$ and allows to exactly compute $t_k$.}

\begin{proof}
Let $k'$ be the route defined by the set of orders $o\in\mathcal{O}$ such that $e_o=1$. 
We need to prove that constraint~\eqref{cst:cutpricing3} is satisfied for $k'$, \ie
$$b  + \sum_{o \in k'} a_o \;+\; \Delta_k \left( |k\cap k'| - |k| + 1 \right) \leq t_{k'}.$$ 
Let us consider the three possible cases:
\begin{enumerate}
\item[case 1]: $k'=k$, then constraint~\eqref{cst:cutpricing3} becomes $t\geq t_k$. It is valid and also permits to exactly compute $t_{k'}$.
\item[case 2]: $k\cap k'=k$, \ie $k\subseteq k'$, then constraint~\eqref{cst:cutpricing3} becomes \\
$t\geq\displaystyle b  + \sum_{o \in k'} a_o \;+\; \Delta_k = t_k+\sum_{o\in k'\setminus k} a_o$\\
This is valid since from hypothesis~\eqref{eq:eq2bhyp}, we obtain $\displaystyle t_{k'}\geq t_k+\sum_{o\in k'\setminus k} a_o.$
\item[case 3]: $k\cap k'\neq k$, then $|k\cap k'| \leq |k| - 1$, so $\Delta_k \left( |k\cap k'| - |k| + 1 \right)\leq 0$. 
However, we have $\displaystyle b  + \sum_{o \in k'} a_o\leq t_{k'}$ from constraint~\eqref{eq:eq2b} applied to route $k'$.
Hence constraint~\eqref{cst:cutpricing3} is valid.
\end{enumerate}
\end{proof}

\subsection{Modifications of the pricing problem due to each type of non-robust cut}
\label{annex_pricing}
\paragraph{Martello and Toth's cuts.}
For each Martello and Toth's cut~\eqref{cnst:l2bound} defined by the pair $(\bar{d}^*, q^*)$, we consider the following additional variables. For the sake of clarity we will not index them by $(\bar{d}^*, q^*)$.

\begin{itemize}
\item $w \geq 0$: real non-negative variable that equal $\bar{d}^*$ if the route is in $\mathcal{K}_1(\bar{d}^*,q^*)$, and equal $t$ if the route is in $\mathcal{K}_2(\bar{d}^*,q^*)$, and 0 otherwise
\item $b_1$: binary variable that equal 1 if the route is in $\mathcal{K}_1(\bar{d}^*,q^*)$, 0 otherwise
\item $b_2$: binary variable that equal 1 if the route is in $\mathcal{K}_2(\bar{d}^*,q^*)$, 0 otherwise
\item $b_3$: binary variable that equal 1 if the route is in $\mathcal{K}_3(\bar{d}^*,q^*)$, 0 otherwise
\end{itemize}

The objective function of $Pr(\mathcal{K}^{'})$ is modified by subtracting $\gamma_{MT}^{\bar{d}^* \, q^*} \, w$, where $\gamma_{MT}^{\bar{d}^* \, q^*}$ is the dual value associated to the cut.
The following constraints are added to $Pr(\mathcal{K}^{'})$:

$$
\begin{array}{cc}
     \displaystyle b_1 + b_2 + b_3 + \sum_{\bar{d} > \bar{d}^*}{\mu^{\bar{d}}} = 1 &
    \constraintlabel{cnst:pri:l2:01}
    \\
    \displaystyle w \geq \bar{d}^* \, b_1 &
    \constraintlabel{cnst:pri:l2:02}
    \\
    \displaystyle w \geq t - \displaystyle\max_{\bar{d} \in \mathcal{D}}{\left\lbrace \bar{d} \right\rbrace} \left( 1 - b_2 \right)&
    \constraintlabel{cnst:pri:l2:03}
    \\
    \displaystyle t \geq q^* \, b_2 + \left( \bar{d}^* - q^* + 1 \right) b_1&
    \constraintlabel{cnst:pri:l2:04} 
    \\
    \displaystyle t \leq \bar{d}^* \, b_1 + \left( \bar{d}^* - q^* \right) b_2 + \left(q^*-1\right) b_3 + \sum_{\bar{d} > \bar{d}^*}{\bar{d} \, \mu^{\bar{d}}}&
    \constraintlabel{cnst:pri:l2:05} 
    \\
    \displaystyle w \geq 0&
    \constraintlabel{cnst:pri:l2:06} 
    \\
    \displaystyle b_1, b_2, b_3 \in \left\lbrace 0, 1 \right\rbrace&
    \constraintlabel{cnst:pri:l2:07} 
\end{array}
$$


Constraint~\eqref{cnst:pri:l2:01} ensures that the route is in one of the three sets $\mathcal{K}_1(\bar{d}^*,q^*)$, $\mathcal{K}_2(\bar{d}^*,q^*)$, $\mathcal{K}_3(\bar{d}^*,q^*)$, or the route has a deadline that is grater than $\bar{d}^*$.
Constraint~\eqref{cnst:pri:l2:02} ensures that $w$ takes value $\bar{d}^*$ if the route is in $\mathcal{K}_1(\bar{d}^*,q^*)$, while constraint~\eqref{cnst:pri:l2:03} ensures that $w$ takes value $t$ if the route is in $\mathcal{K}_2(\bar{d}^*,q^*)$.
Constraints~\eqref{cnst:pri:l2:04} and~\eqref{cnst:pri:l2:05} ensure the total time of the route is consistent with the set that contains the route.
Constraints~\eqref{cnst:pri:l2:06} and~\eqref{cnst:pri:l2:07} define the domain of the decision variables.

\paragraph{Fekete et Schepers cuts.}
For each Fekete and Schepers' cut~\eqref{eq:vi01} defined by the pair $(\bar{d}^*, q^*)$, we need to consider the following additional variable:
\begin{itemize}
\item $v:$ non-negative integer variable that represents the number of intervals of size $\bar{d}^* / q^*$ covered by the total time of the route (with a positive residual), \ie what is denoted by $i$ in constraint~\eqref{eq:vi01}.
\end{itemize}
For the sake of clarity we will not indexed it by $(\bar{d}^*, q^*)$.

The objective function of $Pr(\mathcal{K}^{'})$ is modified by subtracting $\gamma_{FS}^{\bar{d}^* \, q^*} \, v$ where $\gamma_{FS}^{\bar{d}^* \, q^*}$ is the dual value associated to the cut.
The following constraints are added to $Pr(\mathcal{K}^{'})$:

$$
\begin{array}{cc}
\displaystyle t \leq (v+1) \frac{\bar{d}^*}{q^*} + \sum_{\bar{d} > \bar{d}^*}{\left( \bar{d} - \frac{\bar{d}^*}{q^*} \right) \mu^{\bar{d}}}&
\constraintlabel{cnst:pri:fs:01}
\\
\displaystyle t \geq \frac{\bar{d}^*}{q^*} v + \frac{1}{q^*}&
\constraintlabel{cnst:pri:fs:02}
\\
\displaystyle v \in \mathbb{N}&
\constraintlabel{cnst:pri:fs:03} 
\end{array}
$$
Constraints~\eqref{cnst:pri:fs:01} and~\eqref{cnst:pri:fs:02} ensure that $v$ corresponds to the number of intervals of size $\bar{d}^* / q^*$ covered by the travel time of the route, if the deadline of the route is less or equal than $\bar{d}^*$.
Note that the value $1 / q^*$ in constraint~\eqref{cnst:pri:fs:02} is set to ensure $t_k > i \, \bar{d}^* / q^*$ in the definition of $\mathcal{K}(\bar{d}^*,i,q^*)$.
Finally, constraint~\eqref{cnst:pri:fs:03} defines the domain of the variable.

\paragraph{Strengthened capacity cuts.}
For each SCC defined by the set of orders $\mathcal{R}^*$, we need to consider the following additional variable:
\begin{itemize}
    \item $z:$ binary variable that equal 1 if the route contains at least one order in $\mathcal{R}^*$, 0 otherwise.
\end{itemize}
Fort he sake of clarity we do not index the variable by $\mathcal{R}^*$.

The objective function of $Pr(\mathcal{K}')$ is modified by subtracting $\gamma_{SC}^{\mathcal{R}^*} \, z$ where $\gamma_{SC}^{\mathcal{R}^*}$ is the dual value associated to the cut.
The following constraints are added to $Pr(\mathcal{K}')$:
\begin{equation*}
\sum_{o \in \mathcal{R}^*}{e_o} \geq z
\qquad
\constraintlabel{cnst:pri:scc:01}
\end{equation*}

\paragraph{Rank-1 cuts.}
The management of R1Cs in a MIP formulation of the pricing problem has been recently proposed by \cite{Hintsch2021}.
They propose two sets of constraints to be included in the MIP, where each set does not dominate the other and are thus worth consideration.
We detail hereafter these two sets of constraints.

For each R1C cut~\eqref{eq:sr} defined by the pair $(\mathcal{R}^*, \textbf{p}^*)$, we need to consider the following additional variable: $u \in \mathbb{N}$ that represents the coefficient $\left\lfloor \sum_{o \in \mathcal{R}^*}{p_o^* e_o} \right\rfloor$ of the route computed by the pricing problem.
For the sake of clarity we do not index the variable by $(\mathcal{R}^*, \textbf{p}^*)$.

The objective function of $Pr(\mathcal{K}')$ is modified by subtracting $\gamma_{R1}^{\mathcal{R}^* \, \textbf{p}^*} u$ where $\gamma_{R1}^{\mathcal{R}^* \, \textbf{p}^*}$ is the dual value associated to the cut.

Let us consider that the multipliers $\textbf{p}^* = \left\lbrace p^*_1,p^*_2, \ldots, p_q^* \right\rbrace$ can be written as $\left\lbrace \frac{s^*_1}{t^*}, \frac{s^*_2}{t^*}, \ldots, \frac{s_q^*}{t^*} \right\rbrace$ with $s^*_1, s_2^*, \ldots, s_q^*, t^* \in \mathbb{N}^*$.

The first constraint to be added to $Pr(\mathcal{K}')$ is the following:
\begin{equation*}
\sum_{o \in \mathcal{R}^*}{s_o^* e_o} - t^* u \leq t^*-1
\qquad
\constraintlabel{cnst:pri:r1:01}
\end{equation*}

For the second set of constraints, let us introduce the notion of \emph{minimal subset} defined as follows:
a set $\mathcal{M} \subseteq \mathcal{R}^*$ is  a minimal subset for $\mathcal{R}^*$ and multipliers $\textbf{p}^*$ if there exits an integer $m \geq 1$ such that:
\begin{equation*}
\sum_{o \in \mathcal{M}}{p_o} \geq m \qquad \text{ and } \qquad \sum_{o \in \mathcal{M}'}{p_o} < m \qquad \forall \mathcal{M}' \subsetneq \mathcal{M}.
\end{equation*}


Let us denote by $\mathcal{N}^*$ the set of all minimal subsets for $\mathcal{R}^*$ and multipliers $\textbf{p}^*$.
The second set of constraints, one for each element in $\mathcal{N}$, is the following:
\begin{equation*}
\sum_{o \in \mathcal{M}}{e_o} - u \leq |\mathcal{M}| - \left\lfloor \sum_{o \in \mathcal{M}}{p_o^*} \right\rfloor
\qquad \forall \mathcal{M} \in \mathcal{N}^*
\qquad
\constraintlabel{cnst:pri:r1:02}
\end{equation*}

The interested reader is referred to \cite{Hintsch2021} for a detailed explanation, especially on minimal subsets $\mathcal{N}^*$ for each optimal vector of multipliers for R1Cs defined over sets of orders of cardinality between 3 and 5.



\subsection{Effect of using the stabilization from \cite{BriantCCCLO20}}
\label{annex_stabilization}

We conduct experiments to evaluate the effect of using stabilization techniques during column generation.
We compare our version of the \algo{} without stabilization to a version of the \algo{} using stabilization from \cite{BriantCCCLO20} (see section 4.1.4), not modifying their parameter values.
The experiments are carried out on nine instances with 100 orders, a trolley capacity $B=15$, and a warehouse with 12 aisles -- \ie the instances from the first nine rows of Table~\ref{tableLargeRes15}.
All of them can be solved by \algo{} within 4 hours.

Table~\ref{tableStabilization} reports the results.
Each row corresponds to one instance, identified in the first column (\texttt{Id}).
The following columns (\texttt{Cpu\ (s)}) report the computational time of \algo{} without and with the stabilization from \cite{BriantCCCLO20}, respectively.

\begin{table}[ht]
\centering
{\footnotesize
\begin{tabular}{|c||r|r|}
\hline
\multirow{2}{*}{Id} &\multicolumn{2}{|c|}{Cpu (s)} \\
\cline{2-3}
	&	No stabilization	&	With stabilization	\\
\hline
7	&	 3~508 	&	 4~808 	\\
8	&	 7~685 	&	 11~545 	\\
9	&	 4~116 	&	 6~268 	\\
34	&	 3~257 	&	 4~341 	\\
35	&	 1~586 	&	 1~929 	\\
36	&	 3~421 	&	 4~363 	\\
61	&	 9~267 	&	 11~113 	\\
62	&	 2~039 	&	 2~680 	\\
63	&	 6~814 	&	 10~744 	\\
\hline
\end{tabular}
}
\caption{Comparison of computation times with and without stabilization in \algo.}\label{tableStabilization}
\end{table}

Table~\ref{tableStabilization} clearly shows that using stabilization with the parameters defined in \cite{BriantCCCLO20} significantly increases computation times.

\end{document}